# ON THE GEOMETRY OF COMPATIBLE POISSON AND RIEMANNIAN STRUCTURES

NICOLÁS MARTÍNEZ ALBA AND ANDRÉS VARGAS

ABSTRACT. We consider compatibility conditions between Poisson and Riemannian structures on smooth manifolds by means of a contravariant partially complex structure, or $f$-structure, introducing the notion of (almost) Kähler–Poisson manifolds. In addition, we study some of their properties under structure preserving maps and symmetries.

## Contents

## 1. INTRODUCTION

Several notions of compatibility between Poisson and metric structures on manifolds have been proposed in the literature, trying to extend particular properties of standard Kähler manifolds. Associated definitions and terminologies, like holomorphic coordinates, holomorphic Poisson manifolds, foliations by parallel planes, Killing–Poisson and Riemann–Poisson structures, have appeared in the literature, each one generalizing a different viewpoint on Kähler manifolds. Despite the different ways to define useful compatibilities, there is no a unifying theory of manifolds having both Poisson and Riemannian structures. In this paper we propose a notion of compatibility between Poisson, metric, and partially complex structures, jointly with a geometric integrability condition.

Starting from the symplectic structure of a Kähler manifold, several of these generalizations extend this notion to the case of a non-degenerate Poisson structure, via a global metric tensor or defining a local characterization from the metric. Let us mention some of these definitions: Arnild and Huisken [1] introduced a weaker notion of almost Kähler manifold (extended to Poisson manifolds) and studied the curvature tensors, Karabegov [20] used complex Poisson manifolds carrying holomorphic local coordinates, and Mokhov [22] considered compatible Poisson structures of hydrodynamic type. Other sources for generalizations of the notion of an almost complex structure come from the study of foliations and compatible foliated structures as in [13, 19, 27]. In [5], Boucetta introduced a global compatibility condition that involves a metric contravariant connection. If the Poisson tensor is invertible (i.e., defines a symplectic form), and there exists a compatible almost complex structure, this compatibility is precisely equivalent to the existence of a Kähler structure on the manifold.

The basic motivation to introduce a compatibility condition between a Poisson and a metric tensor on a smooth manifold $M$ comes from the geometry behind a symplectic form $\omega$ and a (pseudo) Riemannian metric $\langle -,- \rangle$ on (pseudo) Kähler manifolds. In this case, the compatibility is equivalent to the condition that $\omega$ should be parallel with respect to the Levi-Civita connection $\nabla_\bullet$ associated to $\langle -,- \rangle$, i.e., the covariant tensor field $\nabla_\bullet \omega$ vanishes. This condition can be extended to a compatibility with a Poisson bivector field $\pi$ by mimicking the same idea, namely, to require that $\nabla \pi = 0$ with respect to a suitable connection $\nabla$ chosen according to convenience. Taking $\nabla$ as the above-mentioned covariant Levi-Civita connection $\nabla_\bullet$ determined by $\langle -,- \rangle$, the condition $\nabla_\bullet \pi = 0$ implies that the symplectic foliation of $(M, \pi)$ is regular (see Poisson connections in [29, Thm. 2.20]), but then most of the interesting examples, like linear Poisson structures, are excluded. Indeed, with this kind of condition a regular foliation arises for any covariant torsion-free connection. This forces us to avoid those connections in order to enlarge the family of examples. A good alternative is to consider suitable contravariant connections $\nabla^\bullet \colon \Gamma(E) \to \Gamma(TM) \otimes \Gamma(E)$ compatible with the metric of some Riemannian vector bundle $E \to M$ in an appropriate sense.

If $(M, \pi, \langle -,- \rangle)$ is a Poisson and Riemannian manifold, the Koszul bracket $[-,-]_\pi$ and the dual metric $\langle -,- \rangle$ acting on one-forms determine a unique metric and torsion-free contravariant connection denoted by $\nabla^\bullet$. Due to its definition and properties $\nabla^\bullet$ will be called the contravariant Levi-Civita connection and it is defined through a Koszul-type formula

$$2\langle \nabla^\alpha \beta, \gamma \rangle = \pi^\sharp(\alpha)\langle \beta, \gamma \rangle + \pi^\sharp(\beta)\langle \alpha, \gamma \rangle - \pi^\sharp(\gamma)\langle \alpha, \beta \rangle \\ + \langle [\alpha,\beta]_\pi, \gamma \rangle - \langle [\alpha,\gamma]_\pi, \beta \rangle - \langle [\beta,\gamma]_\pi, \alpha \rangle, \qquad (1.1)$$

where $\pi^\sharp \colon T^*M \to TM$ denotes the anchor map induced by $\pi$. This connection is in fact adequate to define the desired geometric compatibility in the Poisson case. Moreover, once we introduce a suitable bundle map $\mathcal{J} \colon T^*M \to T^*M$ that behaves like a contravariant almost complex structure, or $f$-structure [31], on the cotangent distribution to the leaves of the symplectic foliation, we can impose a contravariantly parallel condition on $\mathcal{J}$ analogous to the covariantly parallel property of a standard complex structure in a Kähler manifold. These are the ingredients needed to introduce the notion of *Kähler–Poisson* structure, that is the central concept in this paper.

One of the main consequences obtained from the contravariantly parallel assumption on $\pi$ is that each leaf of the symplectic foliation inherits a Kähler structure. Moreover, in the regular case the existence of a Kähler foliation turns out to be equivalent to the parallel assumption on $\mathcal{J}$. Old and new results on the geometry of the foliation in this setting will be presented, as well as a detailed study of the behavior of this setup under the presence of symmetries, includying its relation with the corresponding reduced structures. The case of a generalized (almost) Kähler–Poisson structure carrying a non-regular symplectic foliation comes as a more general situation that will be considered in a future work.

The organization of the paper is the following: in Section 2 we begin with a brief review of metric and Poisson structures on manifolds and continue in Section 3 with the introduction of the main notions of the paper, the Kähler–Poisson and Riemann–Poisson structures, together with some of their attributes; in particular we also recover another compatibility, the so-called Killing–Poisson structures (introduced by Boucetta in [5, 7]). In addition to already known results and a characterization of Kähler–Poisson structures in Proposition 3.7, we provide metric consequences for the geometry of the symplectic foliation in Theorems 3.24 and 3.28. In Section 4 we focus on Riemannian submersions that are also Poisson maps, in particular we obtain Propositions 4.5 and Theorem 4.8 that describe the behavior of Riemann–Poisson, Kähler–Poisson, and Killing–Poisson structures under such structure-preserving maps. Finally, in Section 5 we consider the case of symmetries of the Poisson and metric structures under the action of a Lie group, includying conditions to obtain a well-behaved reduced structure. In particular, we state a reduction result in Theorem 5.2, and discuss the case of Kählerian orbits on compact Lie groups.

*Acknowledgments.* We would like to express our gratitude to Henrique Burzstyn for his interest and many insightful comments on preliminary versions of this manuscript that helped to improve our presentation and results. We are also thankful to Rui L. Fernandes and Iván Contreras for comments and suggestions on several topics related to this work. The authors are indebted to the Instituto de Matemática Pura y Aplicada (IMPA), Rio de Janeiro, Brazil, for the hospitality and financial support for visits during the early stages of this project. This work was partially financed by the research project ID-PRY: 6918 of the Faculty of Sciences, Pontificia Universidad Javeriana, Bogotá, Colombia.

## 2. Preliminaries

We begin with a brief summary of Poisson and pseudo-Riemannian structures on smooth manifolds in order to fix our notation and present the basic ingredients needed for the following sections.

2.1. **Poisson structures.** A Poisson manifold is a pair $(M, \pi)$ where $M$ is a smooth manifold and $\pi$ is a bivector field so that $[\pi, \pi]_s = 0$, where $[\text{-},\text{-}]_s$ denotes the Schouten bracket of multivector fields. The Poisson structure is equivalently defined via a Lie algebra $(C^\infty(M), \{\text{-},\text{-}\})$, where the skew-symmetric bracket $\{\text{-},\text{-}\}$ satisfies in addition to the Jacobi identity a Leibniz rule with respect to multiplication of functions, i.e. for all $f, g, h \in C^\infty(M)$ the following relation holds:

$$\{f, gh\} = \{f, g\}h + g\{f, h\}.$$

The bivector $\pi$ and the bracket $\{\text{-},\text{-}\}$ are mutually uniquely determined by the relation

$$\pi(df, dg) = \{f, g\}$$

for any two smooth functions $f, g$ on $M$. In this case the condition $[\pi, \pi]_s = 0$ translates to the Jacobi identity for $\{\text{-},\text{-}\}$, indeed we get that

$$\frac{1}{2}[\pi, \pi]_s(df, dg, dh) = \{\{f, g\}, h\} + \circlearrowleft .$$

The bivector $\pi \in \mathfrak{X}^2(M)$ induces an anchor map $\pi^\sharp \colon T^*M \to TM$, given by $\pi^\sharp(\alpha) := \pi(\alpha, \text{-})$, i.e., $\beta(\pi^\sharp(\alpha)) = \pi(\alpha, \beta)$, for all one-forms $\alpha, \beta \in \Gamma(T^*M)$. For simplicity, where no confusion arises, we will denote $\pi^\sharp \alpha := \pi^\sharp(\alpha)$.

Given a smooth function $t : P \to M$ between two Poisson manifolds, $t$ is called a Poisson map if it preserves the Poisson bracket, i.e., for all $f, g \in C^\infty(M)$

$$\{t^*f, t^*g\}_P = t^*\{f, g\}_M.$$

In terms of the bivectors, the condition of a Poisson map translates to the condition that the vector fields $\pi_P^\sharp(t^*\alpha)$ and $\pi_M^\sharp(\alpha)$ are $t$-related for all one-forms $\alpha$ on $M$.

For any function $f \in C^\infty(M)$ there is an associated vector field $X_f := \{f,\text{-}\} = \pi^\sharp(df)$ called the Hamiltonian vector field associated to $f$. Functions whose Hamiltonian vector fields vanish identically are known as Casimir functions.

The pointwise space spanned by Hamiltonian vector fields coincides with the image of $\pi_p^\sharp \colon T_p^*M \to T_pM$ at each $p \in M$. In addition, the relation $X_{\{f,g\}} = [X_f, X_g]$ implies that the smooth distribution $\operatorname{Im} \pi^\sharp$ is involutive, thus it defines an immersed (possibly singular) foliation $\mathcal{F}$ characterized by the fact that each of its leaves is a connected component of the equivalence relation on points of $M$ given by the existence of joining paths obtained through finite composition of Hamiltonian flows. In addition, each leaf is equipped with a symplectic form defined by $\omega(\pi^\sharp\alpha, \pi^\sharp\beta) := \pi(\alpha,\beta)$. This foliation is known as the *symplectic foliation* of $(M,\pi)$.

We recall also that the Koszul bracket of one-forms $[\alpha,\beta]_\pi := \mathcal{L}_{\pi^\sharp\alpha}\beta - \mathcal{L}_{\pi^\sharp\beta}\alpha - d(\pi(\alpha,\beta))$ endows the space $(\Omega^1(M), [\text{-},\text{-}]_\pi)$ with the structure of a Lie algebra, which will be required in Proposition 4.5.

We conclude this section with two simple lemmata needed to simplify calculations and proofs in the next section.

**Lemma 2.1.** *Let $X \in \Gamma(TM)$ and $\alpha \in \Gamma(T^*M)$. Then*

$$\beta\bigl(\mathcal{L}_X(\pi^\sharp\alpha) - \pi^\sharp(\mathcal{L}_X\alpha)\bigr) = (\mathcal{L}_X\pi)(\alpha,\beta), \tag{2.1}$$

*for any $\beta \in \Gamma(T^*M)$. Therefore, $[\mathcal{L}_X, \pi^\sharp] = 0$ if and only if $\mathcal{L}_X\pi = 0$.*

*Proof.* A straightforward calculation yields

$$\beta\bigl(\mathcal{L}_X(\pi^\sharp\alpha) - \pi^\sharp(\mathcal{L}_X\alpha)\bigr) = \bigl(\mathcal{L}_X(\beta(\pi^\sharp\alpha)) - (\mathcal{L}_X\beta)(\pi^\sharp\alpha)\bigr) - \pi(\mathcal{L}_X\alpha,\beta)$$
$$= \mathcal{L}_X(\pi(\alpha,\beta)) - \pi(\alpha,\mathcal{L}_X\beta) - \pi(\mathcal{L}_X\alpha,\beta) = (\mathcal{L}_X\pi)(\alpha,\beta). \qquad \square$$

Using the definitions of $[\text{-},\text{-}]_\pi$ and $\mathcal{L}_X\pi$, it follows easily from the previous lemma that

**Lemma 2.2.** *Let $X \in \Gamma(TM)$ and $\alpha,\beta \in \Gamma(T^*M)$. Then*

$$[\alpha,\beta]_\pi(X) = \pi^\sharp(\alpha)(\beta(X)) - \pi^\sharp(\beta)(\alpha(X)) + (\mathcal{L}_X\pi)(\alpha,\beta). \tag{2.2}$$

2.2. **Metric structures.** Recall that a pseudo-Riemannian metric tensor $\langle\text{-},\text{-}\rangle$ on $M$, which we will refer simply as a metric structure on $M$, determines two canonical bundle maps called musical isomorphisms: flat $\flat\colon TM \to T^*M$ taking $X \mapsto X^\flat := \langle X,\text{-}\rangle$ and its inverse sharp $\sharp \equiv \flat^{-1}\colon T^*M \to TM$ taking $\alpha \mapsto \alpha^\sharp$, so that $\langle\alpha^\sharp, X\rangle = \alpha(X)$ for any $X \in TM$. In addition, $\langle\text{-},\text{-}\rangle$ has an associated dual metric, or cometric, which will be denoted by slightly curly angle brackets $\langle\!\langle\text{-},\text{-}\rangle\!\rangle$ acting on sections of $T^*M$, that for $X,Y \in \Gamma(TM)$ and $\alpha,\beta \in \Gamma(T^*M)$ satisfies,

$$\langle\!\langle X^\flat, Y^\flat\rangle\!\rangle := \langle X,Y\rangle = X^\flat(Y) \quad \text{and} \quad \langle\!\langle \alpha,\beta\rangle\!\rangle := \langle \alpha^\sharp,\beta^\sharp\rangle = \alpha(\beta^\sharp). \tag{2.3}$$

Furthermore, any pseudo-Riemannian metric $\langle\text{-},\text{-}\rangle$ determines a unique metric and torsion-free covariant derivative $\nabla_\bullet\colon \Gamma(TM) \to \Gamma(T^*M) \otimes \Gamma(TM)$ called the Levi-Civita connection and given explicitly by the Koszul formula

$$2\langle\nabla_X Y, Z\rangle = X\langle Y,Z\rangle + Y\langle X,Z\rangle - Z\langle X,Y\rangle$$
$$+ \langle[X,Y],Z\rangle - \langle[X,Z],Y\rangle - \langle[Y,Z],X\rangle. \tag{2.4}$$

Given a function $f \in C^\infty(M)$ we will denote by $\nabla f$ its gradient with respect to the metric $\langle\text{-},\text{-}\rangle$, determined by the equality $\langle\nabla f, X\rangle = df(X) = Xf$ for all $X \in \mathfrak{X}(M)$. Similarly, given a vector field $X$ we will denote by $\operatorname{div} X \in C^\infty(M)$ its divergence which, assuming $M$ is oriented, can be defined by the condition $\mathcal{L}_X\mu = (\operatorname{div} X)\mu$, where $\mu$ denotes the pseudo-Riemannian volume $n$-form of $M$ written as $\mu := e_1^\flat \wedge \cdots \wedge e_n^\flat$ with respect to a local orthonormal frame $(e_1, \ldots, e_n)$. Equivalently, using the Levi-Civita connection it can be defined as $\operatorname{div} X := \operatorname{tr}(\nabla_\bullet X)$ that has the local expression $\operatorname{div} X = \sum_{i=1}^n \langle \nabla_{e_i} X, e_i \rangle$. If $\iota: N \hookrightarrow M$ is an immersion of a submanifold $N$ of the Riemannian manifold $(M, \langle -, - \rangle)$ we denote the pullback metric on $N$ by $\langle -, - \rangle|_{TN} := \iota^* \langle -, - \rangle$, in particular we will use the notation $\langle -, - \rangle|_{T\mathcal{F}}$ for the metric induced on the leaves of a foliation $\mathcal{F}$ of $(M, \langle -, - \rangle)$.

Finally, general details on contravariant connections can be found in [15, 29]. The Levi-Civita contravariant connection defined by (1.1) is studied in [5, 6, 8].

## 3. Compatibilities

3.1. **Linear case.** Here we will work at the linear algebra level, fixing a finite dimensional real vector space $V$. The aim is to adapt the well known relations among complex, symplectic and inner product structures on $V$ to the degenerate case, i.e. when there is a non-necessarily invertible linear Poisson structure $\pi \in \wedge^2 V$. For this we will fix our setting assuming, according to the case, the existence of the following linear structures:

(p) A *Poisson bivector*: a skew-symmetric bilinear map on forms $\pi \in \wedge^2 V$.
(m) A *dual metric*: a non-degenerate symmetric bilinear map on forms $\langle -, - \rangle \in \operatorname{Sym}^2 V$.
(c) A contravariant *partially complex structure*, or contravariant $f$-*structure*: a linear map $\mathcal{J}: V^* \to V^*$ such that $\mathcal{J}^3 + \mathcal{J} = 0$.

The reader is referred to [31] for details on (covariant) $f$-structures as they were called by K. Yano. The name *hor-complex* structures was also used in [2] but it was not adopted in the literature.

**Definition 3.1.** A triple $(\pi, \langle -, - \rangle, \mathcal{J})$ of linear structures on a vector space $V$ as in (p), (m) and (c), respectively, is called *compatible* if for every $\alpha, \beta \in V^*$ the relation $\pi(\alpha, \beta) = \langle \alpha, \mathcal{J}\beta \rangle$ is satisfied.

Note that the contraction of the dual metric $\langle -, - \rangle$ with a covector corresponds precisely to the linear sharp map $\sharp = \flat^{-1} : V^* \to V$, therefore the linear compatibility condition between $\pi$ and $\mathcal{J}$ can be rephrased as $\pi^\sharp + \sharp \circ \mathcal{J} = 0$ or, more explicitly, $\pi^\sharp(\alpha) = -(\mathcal{J}\alpha)^\sharp$ for all $\alpha \in V^*$.

An analogous notion of compatibility is well-known in the symplectic case and its relation with Definition 3.1 is given in the following

**Proposition 3.2.** *If $\pi$ is invertible so that $\omega := -\pi^{-1}$ is a linear symplectic form, $\langle -, - \rangle$ is an inner product on $V$, and $J: V \to V$ is defined by the relation $\omega(X, Y) = \langle X, JY \rangle$, then $\mathcal{J} \circ \flat \circ J = \flat$. Moreover, $\mathcal{J}^2 = -\operatorname{Id}_{V^*}$ if and only if $J^2 = -\operatorname{Id}_V$.*

In the case of complex, symplectic and metric structures on even dimensional vector spaces it can be proved that any two of these structures uniquely determine the third one (see for example the discussion in [11, Lect. 13]). We state the corresponding result for the case of Poisson structures on linear spaces

**Theorem 3.3.** *On a finite dimensional vector space $V$ the following statements hold:*

(i) *Let $\mathcal{J}$ and $\langle -, - \rangle$ be as in (c) and (m), respectively. If $\mathcal{J} = -\mathcal{J}^t$, then there exists some $\pi$ as in (p) such that the triple is compatible.*
(ii) *Let $\mathcal{J}$ and $\pi$ be as in (c) and (p), respectively. If $\operatorname{Ker} \mathcal{J} = \operatorname{Ker} \pi$ and $\mathcal{J}$ restricts to an isomorphism on the complement of $\operatorname{Ker} \pi$ in $V^*$, then there exists $\langle -, - \rangle$ as in (m) with respect to which $\mathcal{J}$ is skew-symmetric and the triple is compatible.*
(iii) *Let $\langle -, - \rangle$ and $\pi$ be as in (m) and (p), respectively. Then, there exist $\mathcal{J}$ as in (c) so that $\mathcal{J}$ is skew-symmetric w.r.t. some $\langle -, - \rangle_A \in \operatorname{Sym}^2 V$ and the triple $(\pi, \langle -, - \rangle_A, \mathcal{J})$ is compatible.*

*Proof.* The proofs follow the same ideas as in the linear symplectic/complex structure case. For convenience of the reader we present here the argument for property (*iii*).

Let us define $A\colon V^* \to V^*$ via the relation $\pi(\alpha,\beta) = \langle \alpha, A\beta \rangle$ for all $\alpha, \beta \in V^*$. If we set $V_1 := (\operatorname{Ker} \pi^\sharp)^\perp$, we can verify that the restricted map

$$A_1 := A|_{V_1}\colon V_1 \to V_1$$

is an isomorphism. Using a polar decomposition there exists some $J_1\colon V_1 \to V_1$ so that $A_1 = |A_1| J_1$ and $J_1^2 = -\operatorname{Id}_{V_1}$. Finally, on $V^* := \operatorname{Ker} \pi^\sharp \oplus V_1$ we define the partially complex structure $J := 0 \oplus J_1$ and the cometric

$$\langle \alpha, \beta \rangle_A = \begin{cases} \langle \alpha, |A_1|\, \beta \rangle & \text{if } \alpha, \beta \in V_1, \\ \langle \alpha, \beta \rangle & \text{if } \alpha, \beta \in \operatorname{Ker} \pi^\sharp, \\ 0 & \text{in other cases.} \end{cases}$$

with respect to which $J$ is skew symmetric by its definition and it is also easy to verify that the triple $(\pi, \langle \text{-},\text{-} \rangle_A, J)$ is compatible. $\square$

It is important to observe that, as it happens in the linear symplectic case, if we begin with a dual metric $\langle \text{-},\text{-} \rangle$ and a Poisson bivector $\pi$ (as in (m) and (p), respectively) we find a compatible triple $(\pi, \langle \text{-},\text{-} \rangle_A, J)$. But if we apply Theorem 3.3 (*ii*) to $J$ and $\pi$, we do not recover the original dual metric $\langle \text{-},\text{-} \rangle$.

3.2. **Riemann–Poisson and Kähler–Poisson manifolds.** Hereafter, given a smooth manifold $M$ we will denote by $\pi \in \mathfrak{X}^2(M)$ a Poisson bivector, by $\langle \text{-},\text{-} \rangle$ a dual metric on $M$, and by $J\colon T^*M \to T^*M$ a bundle map of fiberwise contravariant partially complex structures, that is $J_p^3 + J_p = 0$ at every point $p \in M$.

Now we consider two different, but related, notions of compatibility. The first one was introduced in [7].

**Definition 3.4.** A manifold $M$ endowed with a pair $(\pi, \langle \text{-},\text{-} \rangle)$ of Poisson and metric structures will be called a *Riemann–Poisson* manifold if $\pi$ is contravariantly parallel, i.e.,

$$\nabla^\bullet \pi = 0, \tag{3.1}$$

where $\nabla^\bullet$ is the contravariant Levi-Civita connection (1.1) associated to $(\pi, \langle \text{-},\text{-} \rangle)$, and $\nabla^\bullet \pi(\alpha, \beta, \gamma) \equiv (\nabla^\alpha \pi)(\beta, \gamma) := \pi^\sharp(\alpha)\pi(\beta, \gamma) - \pi(\nabla^\alpha \beta, \gamma) - \pi(\beta, \nabla^\alpha \gamma)$ for all one-forms $\alpha, \beta, \gamma \in \Gamma(T^*M)$.

**Definition 3.5.** A manifold $M$ equipped with a triple $(\pi, \langle \text{-},\text{-} \rangle, J)$ of Poisson, metric and contravariant partially complex structures will be called an *almost Kähler–Poisson* manifold if at each $p \in M$ the structures $(\pi_p, \langle \text{-},\text{-} \rangle_p, J_p)$ are linearly compatible[1]. In addition, an almost Kähler–Poisson manifold $(M, \pi, \langle \text{-},\text{-} \rangle, J)$ will be called a *Kähler–Poisson* manifold if $J$ is contravariantly parallel, i.e.,

$$\nabla^\bullet J = 0, \tag{3.2}$$

where $\nabla^\bullet J(\alpha, \beta) \equiv (\nabla^\alpha J)(\beta) := \nabla^\alpha(J\beta) - J(\nabla^\alpha \beta)$ for all $\alpha, \beta \in \Gamma(T^*M)$.

Clearly, any Riemannian manifold $(M, g)$ with $\pi = 0$ and $J = 0$ is a trivial example of a Kähler–Poisson manifold with a foliation of 0-dimensional leaves. Another very simple non-trivial example is

---

[1]Henceforth, a Poisson and (pseudo-) Riemannian structure on a manifold $M$ will simply refer to a Poisson bivector field $\pi$ jointly with a (pseudo-) Riemannian metric $\langle \text{-},\text{-} \rangle$ on $M$, while an almost Kähler–Poisson structure will always mean a compatible triple $(\pi, \langle \text{-},\text{-} \rangle, J)$ on $M$.

*Example* 3.6. For any integer $n \geqslant 2$ fixed $r,s \in \{1,\ldots,n\}$ with $r < s$, let us consider $M = \mathbb{R}^n$ endowed with the Poisson structure $\pi_{(rs)} = \partial_r \wedge \partial_s$, and the Euclidean cometric $\langle\text{-},\text{-}\rangle = \sum_{i=1}^n \partial_i \otimes \partial_i$, written with respect to the canonical (global) Cartesian coordinate system $(x_1,\ldots,x_n)$. Since these structures are constant, their Levi-Civita contravariant derivative vanishes. Moreover, they admit the compatible and constant partially complex structure $\mathcal{J}_{(rs)} = \partial_r \otimes dx^s - \partial_s \otimes dx^r$, thus $(M, \pi_{(rs)}, \langle\text{-},\text{-}\rangle, \mathcal{J}_{(rs)})$ is a Kähler–Poisson manifold with a symplectic foliation of 2-dimensional leaves. The case $n = 2$ corresponds to the 2-dimensional Kähler manifold $(\mathbb{R}^2, \pi^{-1}, \langle\text{-},\text{-}\rangle^{-1})$ isometric to $(\mathbb{C}, \omega, h)$ endowed with the Hermitian form $\omega = -\frac{i}{2} dz \wedge d\bar{z}$, and the canonical Hermitian metric $h = \text{Re}\,(dz \otimes d\bar{z})$. The case $n = 3$ corresponds to Euclidean $\mathbb{R}^3$ foliated by parallel $\mathbb{C}$ planes with their own Kähler structure. Moreover, this construction can be generalized to obtain a foliation of $\mathbb{R}^{2k+\ell}$ by $2k$-planes isometric to $\mathbb{C}^k$ with the metric $\text{Re}\,(\sum_{i=1}^k dz_k \otimes d\bar{z}_k)$. $\diamondsuit$

Similarly to the case of an almost complex structure, there is a skew-symmetric contravariant Nijenhuis tensor $N_{\mathcal{J}} \colon \Omega^1(M) \wedge \Omega^1(M) \to \Omega^1(M)$ associated to $\mathcal{J}$ given by

$$N_{\mathcal{J}}(\alpha, \beta) := [\mathcal{J}\alpha, \mathcal{J}\beta]_\pi + \mathcal{J}^2[\alpha,\beta]_\pi - \mathcal{J}\big([\alpha, \mathcal{J}\beta]_\pi + [\mathcal{J}\alpha, \beta]_\pi\big).$$

A characterization of Kähler–Poisson structures and the relation with $N_{\mathcal{J}}$ can now be stated.

**Proposition 3.7.** *An almost Kähler–Poisson manifold $(M, \pi, \langle\text{-},\text{-}\rangle, \mathcal{J})$ is Kähler–Poisson if and only if $\nabla^\bullet \pi = 0$, i.e., for any $\alpha, \beta, \gamma \in \Gamma(T^*M)$ we have*

$$\nabla^\bullet \pi(\alpha,\beta,\gamma) \equiv (\nabla^\alpha \pi)(\beta,\gamma) := \pi^\sharp(\alpha)\pi(\beta,\gamma) - \pi(\nabla^\alpha \beta, \gamma) - \pi(\beta, \nabla^\alpha \gamma) = 0. \qquad (3.3)$$

*Moreover, the contravariant Nijenhuis tensor can be written in terms of the contravariant Levi-Civita connection $\nabla^\bullet$ as*

$$N_{\mathcal{J}}(\alpha, \beta) = \nabla^\bullet \mathcal{J}(\mathcal{J}\alpha, \beta) - \nabla^\bullet \mathcal{J}(\mathcal{J}\beta, \alpha) - \mathcal{J}\big(\nabla^\bullet \mathcal{J}(\alpha,\beta) - \nabla^\bullet \mathcal{J}(\beta,\alpha)\big).$$

*In particular, $\nabla^\bullet \mathcal{J} = 0$ implies that $N_{\mathcal{J}} = 0$.*

Therefore, by this characterization any Kähler–Poisson manifolds is also a Riemann–Poisson manifold.

*Proof.* The first claim follows from the fact that for all $\alpha, \beta, \gamma \in \Gamma(T^*M)$ it holds

$$\begin{aligned}\langle \gamma, (\nabla^\bullet \mathcal{J})(\alpha,\beta)\rangle &= \langle \gamma, \nabla^\alpha(\mathcal{J}\beta)\rangle - \langle \gamma, \mathcal{J}(\nabla^\alpha \beta)\rangle \\ &= \pi^\sharp(\alpha)\langle \gamma, \mathcal{J}\beta\rangle - \langle \nabla^\alpha \gamma, \mathcal{J}\beta\rangle - \langle \gamma, \mathcal{J}(\nabla^\alpha \beta)\rangle \\ &= \pi^\sharp(\alpha)\,\pi(\gamma,\beta) - \pi(\nabla^\alpha \gamma, \beta) - \pi(\gamma, \nabla^\alpha \beta) = (\nabla^\bullet \pi)(\alpha,\beta,\gamma).\end{aligned}$$

Since $\langle\text{-},\text{-}\rangle$ is non-degenerate we conclude that $\nabla^\bullet \mathcal{J} = 0$ if and only if $\nabla^\bullet \pi = 0$.

The second claim on $N_{\mathcal{J}}$ is a straightforward calculation using the torsion-freeness of $\nabla^\bullet$ as it is done with the covariant derivative $\nabla_\bullet$ for standard almost complex structures. $\square$

*Example* 3.8. For $k > 0$, we consider in $\mathbb{R}^3$ the Poisson bivector $\pi = \partial_1 \wedge \partial_2$ and the $k$-rescaled metric $\langle\text{-},\text{-}\rangle_k := k\langle\text{-},\text{-}\rangle$, using the same notation of the previous example. It is routine to verify that the contravariant Levi-Civita connection associated to $(\pi, \langle\text{-},\text{-}\rangle_k)$ coincides with the one associated to $(\pi, \langle\text{-},\text{-}\rangle)$, hence $(\mathbb{R}^3, \pi, \langle\text{-},\text{-}\rangle_k)$ is Riemann–Poisson as a consequence of Proposition 3.7. Finally, note that $(\mathbb{R}^3, \pi, \langle\text{-},\text{-}\rangle, \mathcal{J})$ is Kähler–Poisson while $(\mathbb{R}^3, \pi, \langle\text{-},\text{-}\rangle_k)$ is Rieman–Poisson but not Kähler–Poisson because there exists no partially complex structure $\mathcal{J}_k$ compatible with $\pi$ and $\langle\text{-},\text{-}\rangle_k$. $\diamondsuit$

Now, we study the case of an almost Kähler–Poisson manifold with non-degenerate Poisson bivector field $\pi$. Recall that in this case $\pi$ admits an everywhere well-defined inverse $\omega := -\pi^{-1}$, which is a symplectic form. Our purpose is to explain the relation between a Kähler–Poisson structure and a standard Kähler one. As expected both notions of Levi-Civita parallelism, namely, the covariant $\nabla_\bullet \omega = 0$ and the contravariant $\nabla^\bullet \pi = 0$ ones, coincide in this case. Here we present a sketch of the proof (originally proposed in [5]).

**Proposition 3.9** (Boucetta). *An almost Kähler–Poisson manifold $(M, \pi, \langle -,-\rangle, \mathcal{J})$ with non-degenerate Poisson bivector field $\pi$ is Kähler–Poisson if and only if $(M, \omega, \langle -,-\rangle, J)$ is a Kähler manifold, where $\omega := \pi^{-1}, \langle -,-\rangle := \langle -,-\rangle^{-1}$, and $\mathcal{J} \circ \flat \circ J = \flat$.*

*Proof.* First note that from the relations $\mathcal{J} \circ \flat = -\flat \circ J$ and $\omega^\flat \circ \pi^\sharp = -\operatorname{Id}$ we get that $(\pi, \langle -,-\rangle, \mathcal{J})$ is almost Kähler–Poisson if and only if $(\omega, \langle -,-\rangle, J)$ is almost Kähler, so it is enough to verify that $\nabla^\bullet \pi = 0$ if and only if $\omega$ is parallel with respect to the covariant Levi-Civita connection. For this let us denote by $\nabla^\bullet_\pi$ the contravariant connection obtained from the Levi-Civita covariant connection $\nabla_\bullet$ associated to $\langle -,-\rangle$ by means of $\pi^\sharp$, namely $\nabla^\alpha_\pi \beta := \nabla_{\pi^\sharp \alpha} \beta$.

With this connection it is straightforward to show that $\nabla^\bullet_\pi \pi = 0$ if and only if $\nabla^\bullet \pi = 0$. On the other hand, by direct computation we obtain

$$(\nabla^\bullet_\pi \pi)(\alpha, \beta, \gamma) = -(\nabla_\bullet \omega)(\pi^\sharp \alpha, \pi^\sharp \beta, \pi^\sharp \gamma),$$

for any $\alpha, \beta, \gamma \in \Gamma(T^*M)$, and the claim follows from the fact that $\nabla_\bullet \omega = 0$ for Kähler manifolds and the characterization of Proposition 3.7. $\square$

*Remark* 3.10. The second part of the proof shows that the contravariant parallelism of $\pi$ is in fact equivalent to the covariant parallelism of $\omega$.

We close this section by recalling an interesting fact associated to the existence of a compatible integrable partially complex structure $\mathcal{J}$ for $(M, \pi, \langle -,-\rangle)$. It follows from results presented in [26] that, when the partially complex structure exists, the symplectic foliation has to be regular. Hence, all Kähler–Poisson manifolds have a regular symplectic foliation. In contrast, there exist Riemann–Poisson manifolds with non-regular symplectic foliations (see Examples 3.20, 3.22 and also the conclusion in Section 5.1), thus they are not Kähler–Poisson. In the case of an unique leaf (i.e., the symplectic case), these two compatible conditions are the same:

**Proposition 3.11.** *A symplectic and Riemannian manifold $(M, \omega, \langle -,-\rangle)$ is Kähler if and only if it is Riemann–Poisson.*

*Proof.* From the results in propositions 3.9 and 3.7, it remains to verify that any Riemann–Poisson is also a Kähler manifold. Moreover, the claim in Remark 3.10 says that it is enough to show that if $(M, \omega, \langle -,-\rangle)$ satisfies $\nabla_\bullet \omega = 0$, then there exists a Kähler structure on $M$. First we will denote by $(\omega, \langle -,-\rangle_0, J_0)$ the almost Kähler structure coming from the pair $(\omega, \langle -,-\rangle)$ and recall that $J_0 = A(-A^2)^{-1/2}$ where $A := -\sharp \circ \omega^\flat$. Let $\nabla_\bullet$ be the Levi-Civita covariant connection associated to $\langle -,-\rangle$, thus for any $X, Y, Z$ vector fields on $M$ we have

$$\langle Z, \nabla_X AY \rangle = X \langle Z, AY \rangle - \langle \nabla_X Z, AY \rangle = X\omega(Z, Y) - \omega(\nabla_X Z, Y)$$
$$= \omega(\nabla_X Z, Y) + \omega(Z, \nabla_X Y) - \omega(\nabla_X Z, Y) = \langle Z, A\nabla_X Y \rangle,$$

where the first equality holds because $\nabla_\bullet$ is metric, and the third one comes from the condition $\nabla_\bullet \omega = 0$. We conclude that $A\nabla_X Y = \nabla_X AY$, and from the fact that $J_0$ can be written as polynomial in $A$ [17], we obtain $\nabla_\bullet J_0 = 0$. Now, using that $\omega$ and $J_0$ are parallel with respect to $\nabla_\bullet$ we have

$$X \langle Y, Z \rangle_0 = -X\omega(Y, J_0 Z) = -\omega(\nabla_X Y, J_0 Z) - \omega(Y, \nabla_X (J_0 Z))$$
$$= -\omega(\nabla_X Y, J_0 Z) - \omega(Y, J_0(\nabla_X Z)) = \langle \nabla_X Y, Z \rangle_0 + \langle Y, \nabla_X Z \rangle_0,$$

that is, $\nabla_\bullet$ is torsion-free and metric with respect to $\langle -,-\rangle_0$, thus it is the covariant Levi-Civita connection for $\langle -,-\rangle_0$. Since $(M, \omega, \langle -,-\rangle_0, J_0)$ is almost Kähler and $\nabla_\bullet \omega = 0$, we conclude that it is actually a Kähler manifold. $\square$

3.3. **Properties of the compatible structures.** Let us summarize some useful properties of the contravariant Levi-Civita connection $\nabla^\bullet$, includying as the last one another characterization of Riemann–Poisson structures (cf. [5]–[9]). Short proofs using our assumptions and terminology are presented for convenience.

**Proposition 3.12.** *For a Poisson and Riemannian manifold $(M, \pi, \langle\text{-},\text{-}\rangle)$ and one-forms $\alpha, \beta, \gamma \in \Gamma(T^*M)$ the following properties hold:*

(i) *If $\alpha \in \Gamma(\operatorname{Ker} \pi^\sharp)$ then $\pi^\sharp(\nabla^\alpha \beta) = \pi^\sharp(\nabla^\beta \alpha)$.*

(ii) *If $\nabla^\bullet \pi = 0$ and $\alpha \in \Gamma(\operatorname{Ker} \pi^\sharp)$ then $\pi^\sharp(\nabla^\alpha \bullet) = \pi^\sharp(\nabla^\bullet \alpha) = 0$. Moreover, with these assumptions we actually have $\nabla^\alpha = 0$.*

(iii) $(\mathcal{L}_{\alpha^\sharp} \pi)(\beta, \gamma) = \langle \nabla^\gamma \alpha, \beta \rangle - \langle \nabla^\beta \alpha, \gamma \rangle$.

(iv) $\nabla^\bullet \pi = 0$ *if and only if for all functions $f \in C^\infty(M)$ and one-forms $\beta, \gamma$ we have*

$$\pi(\nabla^\beta df, \gamma) = \pi(\nabla^\gamma df, \beta).$$

*Furthermore, the same results hold for an almost Kähler–Poisson manifold $(M, \pi, \langle\text{-},\text{-}\rangle, \mathcal{J})$ replacing the condition $\nabla^\bullet \pi = 0$ by $\nabla^\bullet \mathcal{J} = 0$.*

*Proof.* Let us denote by $\alpha, \beta, \gamma$ appropriate sections of $T^*M$ as required in each case.

(i) It follows directly from the torsion-freeness of $\nabla^\bullet$.

(ii) From the parallelism we have that for all $\gamma \in \Gamma(T^*M)$ it holds

$$\gamma(\pi^\sharp(\nabla^\beta \alpha)) = \pi^\sharp(\beta)\pi(\alpha, \gamma) - \pi(\alpha, \nabla^\beta \gamma).$$

Hence if $\alpha \in \Gamma(\operatorname{Ker} \pi^\sharp)$ we get that $\pi^\sharp(\nabla^\beta \alpha) = \pi^\sharp(\nabla^\alpha \beta) = 0$ for all $\beta \in \Gamma(T^*M)$ and the first claim holds. For the second one, consider a splitting (defined at least locally on an open subset $U \subset M$):

$$T^*M = \operatorname{Ker} \pi^\sharp \oplus (\operatorname{Ker} \pi^\sharp)^\perp, \tag{3.4}$$

and arbitrary local sections $\beta = \beta_0 \oplus \beta_\perp \in \Gamma_U(T^*M)$ and $\gamma \in \Gamma_U(\operatorname{Ker} \pi^\sharp)$. Since we also have $\nabla^\alpha \gamma \in \Gamma_U(\operatorname{Ker} \pi^\sharp)$, using that $\nabla^\bullet \langle\text{-},\text{-}\rangle = 0$ we find

$$\langle \nabla^\alpha \beta, \gamma \rangle = \pi^\sharp(\alpha)\langle \beta, \gamma \rangle - \langle \beta, \nabla^\alpha \gamma \rangle$$
$$= -\langle \beta_0, \nabla^\alpha \gamma \rangle - \langle \beta_\perp, \nabla^\alpha \gamma \rangle = -\langle \beta_0, \nabla^\alpha \gamma \rangle = 0.$$

where the last equality follows from (1.1) because $\beta_0, \alpha, \gamma \in \Gamma_U(\operatorname{Ker} \pi^\sharp)$. This and the first claim imply that $\nabla^\alpha \beta \in \Gamma_U(\operatorname{Ker} \pi^\sharp \cap (\operatorname{Ker} \pi^\sharp)^\perp)$ which yields the conclusion.

(iii) This is a consequence of the tensoriality of the relation

$$(\mathcal{L}_{\alpha^\sharp} \pi)(df, dg) = \langle \nabla^{dg} \alpha, df \rangle - \langle \nabla^{df} \alpha, dg \rangle,$$

which is proved using the definition of $\mathcal{L}_X \pi$, the fact that $\nabla^\bullet \langle\text{-},\text{-}\rangle = 0$, and the relations

$$\langle df, \alpha \rangle = \mathcal{L}_{\alpha^\sharp} f \quad \text{and} \quad \pi^\sharp(\gamma)\langle \alpha, df \rangle = \pi(\gamma, \mathcal{L}_{\alpha^\sharp} f).$$

(iv) It suffices to verify that for arbitrary functions $f, g, h \in C^\infty(M)$ it holds

$$-\big(\{f, \{g, h\}\} + \circlearrowleft \,\big) = \nabla^\bullet \pi(df, dg, dh) + \pi(\nabla^{dg} df, dh) + \pi(dg, \nabla^{dh} df),$$

and to extend to all one-forms using the tensoriality of $\nabla^\bullet$. □

*Remark* 3.13. Following their proofs, it is clear that the statements of the previous proposition are also true if the almost Kähler–Poisson manifold is replaced by a Poisson and Riemannian manifold $(M, \pi, \langle\text{-},\text{-}\rangle)$, and instead of the condition $\nabla^\bullet \mathcal{J} = 0$ we use the contravariant Poisson connection condition $\nabla^\bullet \pi = 0$.

Note that Proposition 3.12 (*ii*) says that the contravariant parallelism condition for $\pi$ implies that $\nabla^\bullet$ is a closed operation when restricted to the space $\operatorname{Ker}\pi^\sharp$. The same question arises for the complementary space, i.e. under which condition $\nabla^\bullet$ is closed when restricted to act on sections of $(\operatorname{Ker}\pi^\sharp)^\perp$? To provide the answer, consider three one-forms $\alpha, \beta, \gamma \in \Gamma(T^*M)$ with $\alpha, \beta \in (\operatorname{Ker}\pi^\sharp)^\perp$ and $\gamma \in \operatorname{Ker}\pi^\sharp$. Using the defining relation (1.1) for $\nabla^\bullet$ we obtain that

$$2\langle \nabla^\alpha \beta, \gamma \rangle = \langle [\alpha, \beta]_\pi, \gamma \rangle.$$

From this relation we conclude immediately the following

**Theorem 3.14.** *The connection $\nabla^\bullet$ is independent of the transversal dual metric $\langle\text{-},\text{-}\rangle^\perp$ if and only if $(\operatorname{Ker}\pi^\sharp)^\perp$ is involutive with respect to the Koszul bracket $[\text{-},\text{-}]_\pi$.*

*Remark* 3.15. Under the involutivity assumption of the previous proposition, the contravariant connection $\nabla^\bullet$ can be said to be adapted to the foliation, namely, it is an operator

$$\nabla^\bullet \colon \Gamma((\operatorname{Ker}\pi)^\perp)^* \to \Gamma(\operatorname{Ker}\pi)^\perp \otimes \Gamma((\operatorname{Ker}\pi)^\perp)^*$$

and this assumption will also simplify the proofs in Proposition 3.12.

Two interesting properties of Kähler–Poisson manifolds are related to the behavior of the bivector $\pi$. One is its invariance along gradient vector fields arising from Casimir functions, and the other is the vanishing of its divergence. Recall that the divergence $\operatorname{div}\pi$ of $\pi$ is defined as the unique vector field satisfying $(\operatorname{div}\pi)(f) = \operatorname{div}(X_f)$ for all smooth functions $f \in C^\infty(M)$. The following proposition is the Kähler–Poisson version of [6, Thm. 1.3], whose proof is presented here for convenience.

**Proposition 3.16.** *If $(M, \pi, \langle\text{-},\text{-}\rangle, \mathcal{J})$ is Riemann–Poisson or Kähler–Poisson then $\pi$ satisfies the following properties,*

(*i*) *Transversal invariance:* $\mathcal{L}_{\nabla f}\pi = 0$, *for any Casimir function $f \in C^\infty(M)$.*
(*ii*) *Vanishing divergence:* $\operatorname{div}\pi = 0$.

*Proof.* For condition (*i*) let us suppose that $f \in C^\infty(M)$ is a Casimir function, then $\nabla^{df} = 0$ by Proposition 3.12 (*ii*). Moreover, for any local section $\beta$ of $T^*M$ it holds $[df, \beta]_\pi = 0$. Hence $\nabla^\beta df = \nabla^{df}\beta + [\beta, df]_\pi = 0$, and thus $\nabla^\bullet df = 0$. Finally, by applying item (*iv*) in Proposition 3.12 with $\alpha = df$ we conclude that $\mathcal{L}_{\nabla f}\pi = 0$.

Now we sketch the proof of condition (*ii*). For this, we fix a local basis $(\alpha_1, \ldots, \alpha_n)$ of $T^*M$ with respect to the splitting (3.4). From [8, eq. (13)] we have the expression $(\operatorname{div}\pi)(f) = \sum_{i=1}^n \langle \nabla^{\alpha_i}df, \alpha_i \rangle$. By using a local basis $(\alpha_1, \ldots, \alpha_l)$ of $\operatorname{Ker}\pi^\sharp$ and the metric $\langle\text{-},\text{-}\rangle_\mathcal{F} := \langle\text{-},\text{-}\rangle|_{T\mathcal{F}}$ along the leaves of $\mathcal{F}$ we get

$$(\operatorname{div}\pi)(f) = \sum_{i=l+1}^n \langle \nabla^{\alpha_i}df, \alpha_i \rangle = \sum_{i=l+1}^n \langle \nabla_{\pi^\sharp \alpha_i}df, \pi^\sharp \alpha_i \rangle_\mathcal{F} = -\operatorname{div}_\mathcal{F}(X_f).$$

The conclusion follows because the right-hand side vanishes by (3.5) applied to $\nu = \omega^n$ and the fact that for any Hamiltonian vector field the property $\mathcal{L}_{X_f}\omega = 0$ holds. □

Note that both conditions in the previous proposition just require the use of the Poisson and the pseudo-Riemannian structures, therefore, they are independent of the geometrical quantity $\mathcal{J}$. According to [6, Def. 1.1], a Poisson manifold $(M, \pi)$ equipped with a pseudo-Riemannian metric $\langle\text{-},\text{-}\rangle$ so that conditions (*i*) and (*ii*) of the previous proposition hold is called a *Killing–Poisson* manifold. Equivalently, we say that a transversally invariant and divergence-free Poisson bivector $\pi$ defines a *Killing–Poisson* structure on $(M, \langle\text{-},\text{-}\rangle)$. Concerning transversal invariance, we have the additional characterization below.

**Proposition 3.17.** *Let $(M, \pi, \langle \text{-},\text{-}\rangle)$ be a Poisson manifold carrying a (pseudo) Riemannian metric $\langle \text{-},\text{-}\rangle$ and suppose that $\mathcal{F}$ denotes its canonical symplectic foliation. Then the following conditions are equivalent:*

(i) $\mathcal{L}_X \pi = 0$ *for all* $X \in \Gamma(T\mathcal{F}^\perp)$.
(ii) $(\operatorname{Ker} \pi^\sharp)^\perp$ *is involutive with respect to* $[\text{-},\text{-}]_\pi$.

*Proof.* Let $\alpha, \beta, \gamma$ be one-forms so that $\alpha, \beta \in (\operatorname{Ker} \pi^\sharp)^\perp$, and $\gamma \in \operatorname{Ker} \pi^\sharp$. By Lemma 2.2 we have immediately that

$$\langle [\alpha, \beta]_\pi, \gamma \rangle = [\alpha, \beta]_\pi(\gamma^\sharp) = (\mathcal{L}_{\gamma^\sharp} \pi)(\alpha, \beta).$$

The equivalence follows by noticing that if $\gamma \in \operatorname{Ker} \pi^\sharp$ the vector field $X := \gamma^\sharp \in \Gamma(T\mathcal{F}^\perp)$, and conversely, if $X \in \Gamma(T\mathcal{F}^\perp)$ the one-form $\gamma := X^\flat \in \operatorname{Ker} \pi^\sharp$ and $\gamma^\sharp = X$. $\square$

*Remark* 3.18. On an oriented Poisson manifold $(M, \pi)$ with volume form $\nu$ the divergence $\operatorname{div}_\nu \pi$, denoted simply by $\operatorname{div} \pi$ when $\nu$ is clear from the context, corresponds precisely to the negative of the modular vector field (cf. [21, Prop. 4.17]), i.e.,

$$\mathcal{L}_{X_f} \nu = -(\operatorname{div} \pi)(f)\, \nu. \tag{3.5}$$

Therefore, a useful characterization of the vanishing divergence condition for the Poisson bivector is given by the following equivalence:

$$\operatorname{div} \pi = 0 \quad \Longleftrightarrow \quad \mathcal{L}_{X_f} \nu = 0 \ \text{ for all } f \in C^\infty(M).$$

In fact, it is well-known that a Poisson manifold $(M, \pi)$ admits a volume form $\nu$ preserved by all Hamiltonian flows if and only if its modular vector field vanishes, and such a volume form is unique up to multiplication by a nowhere vanishing Casimir function. Volume forms with this property are called *invariant densities* and are only available on unimodular Poisson manifolds (see [21] for details).

*Example* 3.19 (Symplectic case). In the symplectic case, i.e. when $\pi$ is invertible and $\omega := \pi^{-1}$, it is easy to verify that for every metric on the $2n$-dimensional manifold $M$, transversal invariance of Proposition 3.16 (i) is trivially true because the Casimir functions are the constant ones. If $\langle \text{-},\text{-}\rangle$ is a Riemannian metric on $M$ with associated Riemannian volume form $\nu$, there exists a nowhere vanishing function $\rho \in C^\infty(M)$ such that $\nu = \rho\, \omega^n$. From Remark 3.18 it follows that $\operatorname{div} \pi = 0$ if and only if $d\rho = 0$. Hence, if we define the function $\phi := \frac{1}{n} \log |\rho|$, the conformally equivalent metric $e^{2\phi} \langle \text{-},\text{-}\rangle = |\rho|^{\frac{2}{n}} \langle \text{-},\text{-}\rangle$ turns the triple $(M, \pi, e^{2\phi}\langle \text{-},\text{-}\rangle)$ into a Killing–Poisson manifold. Alternatively, using this notation we also have that $(M, |\rho|^{-1/n}\pi, \langle \text{-},\text{-}\rangle)$ is Killing–Poisson. $\diamondsuit$

The previous example reveals an important difference between Killing–Poisson and Kähler–Poisson structures (see Proposition 3.11). Indeed, it shows that the Killing–Poisson condition relies solely on the metric and not on the partially complex structure as the Kähler–Poisson condition does. The remaining examples are illustrated in Figure 3.3.

*Example* 3.20 (Lie–Poisson structure of $\mathfrak{so}_3^*(\mathbb{R})$). Consider the dual Lie algebra $\mathfrak{so}_3^*(\mathbb{R})$, which as a vector space is isomorphic to $\mathbb{R}^3$, with its Lie–Poisson structure defined by the bivector

$$\pi_{\mathfrak{so}_3^*} := z\partial_x \wedge \partial_y - y\partial_z \wedge \partial_x + x\partial_y \wedge \partial_z.$$

The rank of $\pi_{\mathfrak{so}_3^*}$ is two at any point different from the origin $o := (0,0,0) \in \mathbb{R}^3$, where it vanishes. The symplectic leaves passing through a point $p := (x, y, z) \in \mathbb{R}^3$ are given by the two-dimensional spheres $S_r^2 := \{(x, y, z) \in \mathbb{R}^3 \mid x^2 + y^2 + z^2 = r^2\}$ with $r \in \mathbb{R}_+$, except at the origin $o$, where the leaf corresponds to the point when $r = 0$. Therefore, the symplectic foliation is singular at the origin.

If $\langle \text{-},\text{-}\rangle_{\mathbb{R}^3}$ denotes the standard three-dimensional Euclidean metric, a straightforward calculation shows that $\nabla^\bullet \pi_{\mathfrak{so}_3^*} = 0$, so that the triple $(\mathfrak{so}_3^*(\mathbb{R}), \pi_{\mathfrak{so}_3^*}, \langle \text{-},\text{-}\rangle_{\mathbb{R}^3})$ is a Riemann-Poisson manifold. Moreover, the rescaled Poisson bivector $\widetilde{\pi} := r\pi_{\mathfrak{so}_3^*(\mathbb{R})}$, where the function $r \equiv r(x,y,z) := (x^2 + y^2 + z^2)^{\frac{1}{2}}$, turns out to be transversally invariant and divergence free with respect to $\langle \cdot,\cdot \rangle_{\mathbb{R}^3}$, and therefore defines a Killing–Poisson structure on the Riemannian manifold $(\mathfrak{so}_3^*(\mathbb{R}), \langle \cdot,\cdot \rangle_{\mathbb{R}^3})$. ◇

*Example* 3.21 (Regular Poisson structure on Euclidean $\mathbb{R}^3 \setminus \{o\}$). Based on the previous example and using the same notation, we restrict the Lie–Poisson structure $\pi_{\mathfrak{so}_3^*}$ of the dual Lie-algebra $\mathfrak{so}_3^*(\mathbb{R}) \cong \mathbb{R}^3$ to its regular part in order to obtain a regular Poisson structure $\pi_{\mathrm{reg}} := \pi_{\mathfrak{so}_3^*}|_{\mathbb{R}^3 \setminus \{o\}}$ on $\mathbb{R}^3 \setminus \{o\}$. Then, the foliation of $\mathbb{R}^3 \setminus \{o\}$ by concentric 2-spheres is regular and with the conformally Euclidean metric $\langle \cdot,\cdot \rangle_r := r^{-1}\langle \cdot,\cdot \rangle_{\mathbb{R}^3}$ the triple $(\mathbb{R}^3 \setminus \{o\}, \pi_{\mathrm{reg}}, \langle \cdot,\cdot \rangle_r)$ is again a Killing–Poisson manifold. Furthermore, defining $\mathcal{J}_{\mathrm{reg}} := z\partial_x \wedge dy - y\partial_z \wedge dx + x\partial_y \wedge dz$, the structure $(\mathbb{R}^3 \setminus \{o\}, \pi_{\mathrm{reg}}, \mathcal{J}_{\mathrm{reg}}, \langle \cdot,\cdot \rangle_r)$ is a Kähler–Poisson manifold, where $\langle \cdot,\cdot \rangle_r$ denotes the dual metric associated to $\langle \cdot,\cdot \rangle_r$. ◇

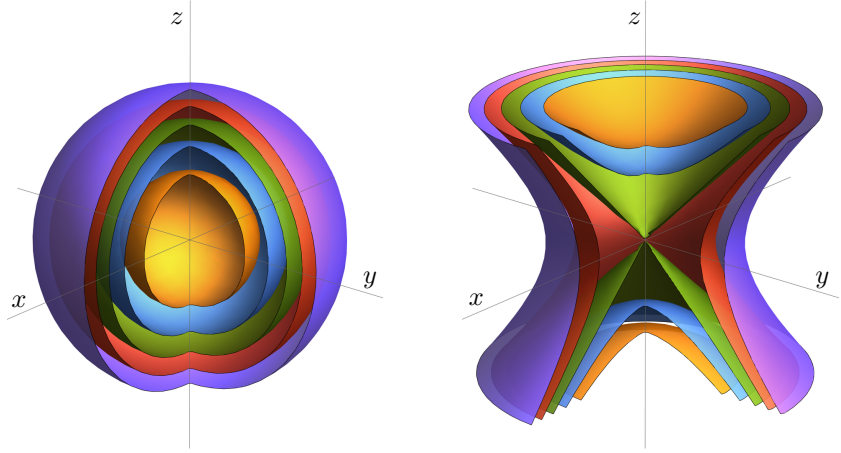

FIGURE 1. Symplectic foliations of Euclidean $\mathbb{R}^3$ and Lorentzian $\mathbb{R}^{2,1}$.
(a) Left: foliation of $(\mathfrak{so}_3^*(\mathbb{R}), \pi_{\mathfrak{so}_3^*(\mathbb{R})}, \langle \cdot,\cdot \rangle_{\mathbb{R}^3})$ by concentric spheres.
(b) Right: foliation of $(\mathfrak{sl}_2^*(\mathbb{R}), \pi_{\mathfrak{sl}_2^*(\mathbb{R})}, \langle \cdot,\cdot \rangle_{\mathbb{R}^{2,1}})$ by coaxial hyperboloids.

*Example* 3.22 (Lie–Poisson structure of $\mathfrak{sl}_2^*(\mathbb{R})$). Now, let us consider the dual Lie algebra $\mathfrak{sl}_2^*(\mathbb{R})$, which is also isomorphic as a vector space to $\mathbb{R}^3$, with its Lie–Poisson structure

$$\pi_{\mathfrak{sl}_2^*} := -z\partial_x \wedge \partial_y + y\partial_z \wedge \partial_x + x\partial_y \wedge \partial_z.$$

The rank of $\pi_{\mathfrak{sl}_2^*}$ is two at any point different from the origin $o := (0,0,0) \in \mathbb{R}^3$, where it vanishes. The symplectic leaves passing through a point $p := (x,y,z) \in \mathbb{R}^3$ are given by the two-dimensional hyperboloids $H_{\pm\rho}^2 := \{(x,y,z) \in \mathbb{R}^3 \mid x^2 + y^2 - z^2 = \pm\rho^2\}$ with $\rho \geqslant 0$ (which are actually cones for $\rho = 0$), except at the origin $o$ that is itself a 0-dimensional leaf (the vertex of the cone $\rho = 0$). Now, denoting by $\mathbb{R}^{2,1}$ the vector space $\mathbb{R}^3$ endowed with the three-dimensional Lorentzian metric $\langle \cdot,\cdot \rangle_{\mathbb{R}^{2,1}} := dx^2 + dy^2 - dz^2$, it follows that $\nabla^\bullet \pi_{\mathfrak{sl}_2^*} = 0$, so that the triple $(\mathfrak{sl}_2^*(\mathbb{R}), \pi_{\mathfrak{sl}_2^*}, \langle \cdot,\cdot \rangle_{\mathbb{R}^{2,1}})$ is a pseudo-Riemann–Poisson manifold. Moreover, the rescaled Poisson bivector $\widetilde{\pi} := \rho\,\pi_{\mathfrak{sl}_2^*(\mathbb{R})}$, where $\rho \equiv \rho(x,y,z) := (x^2 + y^2 - z^2)^{\frac{1}{2}}$, defines a Killing–Poisson structure on the Lorentzian manifold $(\mathfrak{sl}_2^*(\mathbb{R}), \langle \cdot,\cdot \rangle_{\mathbb{R}^{2,1}})$. ◇

*Example* 3.23 (Regular Poisson structure on Lorentzian $\mathbb{R}^{2,1} \setminus H_0^2$). If we restrict the Lie–Poisson structure $\pi_{\mathfrak{sl}_2^*}$ of $\mathfrak{sl}_2^*(\mathbb{R})$ to the regular part lying outside of the cone $H_0^2$, we obtain a regular Poisson structure $\pi_{\mathrm{reg}}$ on $\mathbb{R}^{2,1} \setminus H_0^2$. The foliation by coaxial hyperboloids is regular and with the conformally equivalent metric $\langle\text{-},\text{-}\rangle_\rho := \rho^{-1}\langle\text{-},\text{-}\rangle_{\mathbb{R}^{2,1}}$ restricted to $\mathbb{R}^{2,1}\setminus H_0^2$ the bivector $\pi_{\text{reg}}$ turns out to be a Killing–Poisson structure. Furthermore, defining $\mathcal{J}_{\text{reg}} := -z\partial_x\wedge dy+y\partial_z\wedge dx-x\partial_y\wedge dz$, the structure $(\mathbb{R}^{2,1}\setminus H_0^2, \pi_{\text{reg}}, \mathcal{J}_{\text{reg}}, \langle\text{-},\text{-}\rangle_\rho)$ is a pseudo-Kähler–Poisson manifold, where $\langle\text{-},\text{-}\rangle_\rho$ denotes the dual metric associated to $\langle\text{-},\text{-}\rangle_\rho$.

$\diamondsuit$

3.4. **Geometry of the symplectic foliation.** Before stating the main results of this section, let us briefly review some basic notions from the theory of (pseudo-Riemannian) foliations and introduce some notation. Given a foliation $\mathcal{F}$ for a smooth manifold $M$, a *transverse pseudo-Riemannian metric* on $(M,\mathcal{F})$ is a symmetric 2-covariant tensor field denoted by $\langle\text{-},\text{-}\rangle^\perp\colon \mathfrak{X}(M)\times\mathfrak{X}(M)\to C^\infty(M)$ that satisfies the following conditions:

(i) $\text{Ker}\,(\langle\text{-},\text{-}\rangle_p^\perp) = T_p\mathcal{F}$, for any $p\in M$, and
(ii) $\mathcal{L}_X\langle\text{-},\text{-}\rangle^\perp = 0$, for any $X\in\Gamma(T\mathcal{F})$.

Property (ii) is known as *(infinitesimal) holonomy invariance* of $\langle\text{-},\text{-}\rangle^\perp$ and in this case the triple $(M,\langle\text{-},\text{-}\rangle,\mathcal{F})$ is called a pseudo-Riemannian foliation. The Riemannian case only requires the extra-assumption that $\langle\text{-},\text{-}\rangle$ is positive definite, i.e. $\langle X,X\rangle\geqslant 0$ for all $X\in\mathfrak{X}(M)$.

Given a Riemannian metric $\langle\text{-},\text{-}\rangle$ on $M$ there is a canonical direct sum decomposition $TM = T\mathcal{F}\oplus T\mathcal{F}^\perp$, where $T\mathcal{F}$ denotes the subbundle of $TM$ tangent to the leaves of the foliation $\mathcal{F}$ and $T\mathcal{F}^\perp := (T\mathcal{F})^\perp$ denotes its orthogonal subbundle. Moreover, the metric inherits a corresponding decomposition

$$\langle\text{-},\text{-}\rangle = \langle\text{-},\text{-}\rangle^\top + \langle\text{-},\text{-}\rangle^\perp, \tag{3.6}$$

where each term is a Riemannian metric when restricted to the subbundles $T\mathcal{F}$ and $T\mathcal{F}^\perp$, respectively. Let $N\mathcal{F} := TM/T\mathcal{F}$ denote the normal subbundle of $\mathcal{F}$ in $M$ arising from the exact sequence of vector bundles

$$0\to T\mathcal{F}\xrightarrow{\iota} TM \xrightarrow{\text{pr}^\perp} N\mathcal{F}\to 0.$$

Here, the notation $\text{pr}^\perp$ for the projection onto the normal subbundle is justified by the isomorphism $N\mathcal{F}\cong T\mathcal{F}^\perp$.

In the pseudo-Riemannian context the explicit assumption

$$T\mathcal{F}\cap T\mathcal{F}^\perp = \{0\} \tag{3.7}$$

is needed in order to obtain an analogous direct sum decomposition of $TM$. This condition will always be assumed throughout the text, making the pseudo-Riemannian setting look formally like the Riemannian one. Finally, a pseudo-Riemannian metric $\langle\text{-},\text{-}\rangle$ is called a *bundle-like* metric with respect to $\mathcal{F}$ if the orthogonal part $\langle\text{-},\text{-}\rangle^\perp$ of its decomposition is a transverse metric.

Now, we present a geometric property of symplectic foliations of Poisson manifolds in the presence of metrics with appropriate compatibility.

**Theorem 3.24.** *Let $(M,\pi,\langle\text{-},\text{-}\rangle)$ be a Poisson and Riemannian manifold. Then:*

(i) *If the Poisson bivector $\pi$ is transversally invariant then $\langle\text{-},\text{-}\rangle$ is a bundle-like metric, and endowed with the transverse metric $\langle\text{-},\text{-}\rangle^\perp$ the symplectic foliation of $(M,\pi,\langle\text{-},\text{-}\rangle^\perp)$ is a Riemannian foliation.*
(ii) *If $(M,\langle\text{-},\text{-}\rangle,\pi)$ is Riemann–Poisson, each leaf admits a Kähler metric compatible with its symplectic structure.*
(iii) *If in addition, $(M,\langle\text{-},\text{-}\rangle,\pi,\mathcal{J})$ is Kähler–Poisson then the Kähler metric on the leaves is precisely the restricted one from $M$.*

A straightforward corollary is that Kähler–Poisson manifolds carry Kählerian foliations, which explains the terminology.

*Proof.* (*i*) Let $\langle -,- \rangle^\perp$ denote the orthogonal part of the Riemannian metric $\langle -,- \rangle$ in the decomposition (3.6). Note that by definition $\langle -,- \rangle^\perp$ only vanishes identically along $T\mathcal{F}$, the tangent distribution to the foliation. Therefore, to prove that $\langle -,- \rangle$ is bundle-like it only remains to show that $\mathcal{L}_X \langle -,- \rangle^\perp = 0$ for all Hamiltonian vector fields $X \in \mathfrak{X}(M)$. Taking two Casimir functions $u, v \in C^\infty(M)$, so that $X_f(u) = X_f(v) = 0$ for any $f \in C^\infty(M)$, we calculate

$$(\mathcal{L}_{X_f}\langle -,- \rangle^\perp)(\nabla u, \nabla v) = X_f \langle \nabla u, \nabla v \rangle - \langle [X_f, \nabla u], \nabla v \rangle - \langle [X_f, \nabla v], \nabla u \rangle$$
$$= X_f(\nabla u(v)) - [X_f, \nabla u](v) - [X_f, \nabla v](u)$$
$$= \nabla u(X_f(v)) - X_f(\nabla v(u)) + \nabla v(X_f(u)) = -X_f(\nabla v(u)) = -\{f, \nabla v(u)\}$$
$$= -\nabla v(\{f, u\}) + \{\nabla v(f), u\} = \nabla v(X_f(u)) - \{X_f(v), u\} = 0.$$

Thus, by definition, the symplectic foliation of $(M, \pi, \langle -,- \rangle^\perp)$ is Riemannian.

(*ii*) For this claim we consider the splitting of $T^*M$ in (3.4), and fixing a leaf of $\mathcal{F}$ we get the isomorphism $\pi_\mathcal{F}^\sharp \equiv \pi^\sharp|_{(\text{Ker } \pi^\sharp)^\perp} : (\text{Ker } \pi^\sharp)^\perp \to T\mathcal{F}$. We now take as a metric along tangent directions to $\mathcal{F}$ the one defined for $X, Y \in \Gamma(T\mathcal{F})$ by

$$\langle X, Y \rangle_\mathcal{F} := \langle (\pi_\mathcal{F}^\sharp)^{-1}(X), (\pi_\mathcal{F}^\sharp)^{-1}(Y) \rangle |_{(\text{Ker } \pi^\sharp)^\perp}. \tag{3.8}$$

The covariant Levi-Civita connection of $\langle -,- \rangle_\mathcal{F}$ is given, on vector fields, by

$$\nabla_\bullet^\mathcal{F} : \Gamma(T\mathcal{F}_p) \times \Gamma(T\mathcal{F}_p) \longrightarrow \Gamma(T\mathcal{F})$$
$$(\pi^\sharp(\alpha), \pi^\sharp(\beta)) \longmapsto \nabla_{\pi^\sharp \alpha}^\mathcal{F} \pi^\sharp \beta := \pi^\sharp(\nabla^\alpha \beta).$$

From a straightforward calculation it follows that the symplectic form $\omega_\mathcal{F} = -\pi_\mathcal{F}^{-1}$ on each leaf of $\mathcal{F}$ is parallel with respect to $\nabla_\bullet^\mathcal{F}$. Here we employ the same argument as in Proposition 3.11 to conclude that there exist a metric $\langle -,- \rangle_0$ and a complex structure $J_0$ that turn $(\mathcal{F}, \omega_\mathcal{F}, \langle -,- \rangle_0, J_0)$ into a Kähler manifold.

(*iii*) Since $(M, \pi, \langle -,- \rangle, \mathcal{J})$ is assumed to be a Kähler–Poisson manifold we can make use of the previous claim first, and to prove this statement it only remains to verify that $\langle -,- \rangle_\mathcal{F}$ is indeed $\langle -,- \rangle|_{T\mathcal{F}}$, the restriction of $\langle -,- \rangle$ along each leaf of $\mathcal{F}$. For this, we invoke Lemma 3.25 below. $\square$

**Lemma 3.25.** *The structures defining an almost Kähler–Poisson manifold $(M, \pi, \langle -,- \rangle, \mathcal{J})$ satisfy the following properties:*

  (*i*) $\mathcal{J}$ *restricts to $T^*\mathcal{F}$. Furthermore, if $\flat_\mathcal{F} : T\mathcal{F} \to T^*\mathcal{F}$ denotes the isomorphism determined by the restricted metric $\langle -,- \rangle|_{T\mathcal{F}}$ and $\sharp_\mathcal{F} = \flat_\mathcal{F}^{-1}$, then the composition $J := -\sharp_\mathcal{F} \circ \mathcal{J} \circ \flat_\mathcal{F} : T\mathcal{F} \to T\mathcal{F}$ is an almost complex structure on the leaves of $\mathcal{F}$.*
  (*ii*) *If $\langle -,- \rangle|_{T\mathcal{F}}$ denotes the restriction of the metric to sections of $T\mathcal{F}$, then*

$$\omega_\mathcal{F}(X, Y) = \langle X, JY \rangle|_{T\mathcal{F}}.$$

  (*iii*) *The metric in (3.8) coincides with $\langle -,- \rangle|_{T\mathcal{F}}$.*

*Proof.* (*i*) This is a straightforward verification using, as in Theorem 3.3 (*ii*), the fact that $\text{Ker } \mathcal{J} = \text{Ker } \pi^\sharp$.

(*ii*) Taking $Y := \pi^\sharp(\beta) = -(\mathcal{J}\beta)^\sharp$ for some $\beta \in (\text{Ker } \pi^\sharp)^\perp$, we have that $\omega_\mathcal{F}(X, Y) = i_X \beta$. Now, the result follows from the equalities

$$i_X \beta = -i_X \mathcal{J}^2 \beta = -i_X \mathcal{J}((\mathcal{J}\beta)^{\sharp_\mathcal{F}})^{\flat_\mathcal{F}} = i_X \mathcal{J}(Y^\flat) = \langle X, (\mathcal{J}(Y^{\flat_\mathcal{F}}))^{\sharp_\mathcal{F}} \rangle|_{T\mathcal{F}} = \langle X, JY \rangle|_{T\mathcal{F}}.$$

(*iii*) A simple calculation yields

$$\langle X, Y \rangle_\mathcal{F} = \langle \omega^\flat X, \omega^\flat Y \rangle = i_{(\omega^\flat X)^{\sharp_\mathcal{F}}} \omega^\flat Y = \omega(Y, (\omega^\flat X)^{\sharp_\mathcal{F}})$$
$$= -\omega((\omega^\flat X)^{\sharp_\mathcal{F}}, Y) = -\langle (\omega^\flat X)^{\sharp_\mathcal{F}}, JY \rangle|_{T\mathcal{F}} = -\omega(X, JY)$$
$$= -\langle X, J^2 Y \rangle|_{T\mathcal{F}} = \langle X, Y \rangle|_{T\mathcal{F}}. \qquad \square$$

*Remark* 3.26. Let us stress the difference between a Riemann–Poisson structure (a Poisson and Riemannian manifold $(M, \pi, \langle\text{-},\text{-}\rangle)$ so that $\nabla^\bullet \pi = 0$) and a Kähler–Poisson one. Both structures induce Kähler metrics on the leaves but only for the Kähler–Poisson case those metrics are of the form $\langle\text{-},\text{-}\rangle|_{T\mathcal{F}}$, i.e., obtained by the restriction of the global metric $\langle\text{-},\text{-}\rangle$ acting on $TM$ to the subbundle $T\mathcal{F}$.

In the case of Poisson manifolds with regular foliations we have a one-to-one correspondence between Kähler–Poisson structures and the admissible Kählerian structures on the leaves as a consequence of the following result

**Proposition 3.27.** *Let $(M, \pi)$ be a regular Poisson manifold such that its canonical symplectic foliation is a Riemannian foliation with Kählerian leaves. Then there exists a dual metric $\langle\text{-},\text{-}\rangle$ on $M$ such that the structure $(M, \pi, \langle\text{-},\text{-}\rangle, \mathcal{J})$ is Kähler–Poisson.*

*Proof.* Let $\mathcal{F}$ denote the symplectic foliation, $T\mathcal{F} = \pi^\sharp(T^*M)$ its tangent distribution, $\langle\text{-},\text{-}\rangle_{\mathcal{F}}$ the Kähler metric on its leaves and $\langle\text{-},\text{-}\rangle^\perp$ the transverse metric of the Riemannian foliation. Since $\mathcal{F}$ is regular, we can construct a proper metric

$$\langle\text{-},\text{-}\rangle_M = \langle\text{-},\text{-}\rangle_{\mathcal{F}} \oplus \langle\text{-},\text{-}\rangle^\perp|_{N\mathcal{F}}$$

defined on the whole of $M$ and obtain well-defined isomorphisms $\operatorname{Ker}\pi^\sharp \cong \operatorname{Ann}(T\mathcal{F})$ for the annihilator of $T\mathcal{F}$, and $\pi^\sharp|_{(\operatorname{Ker}\pi^\sharp)^\perp} \colon (\operatorname{Ker}\pi^\sharp)^\perp \to T\mathcal{F}$, where $\perp$ refers to the notion of orthogonality determined by $\langle\text{-},\text{-}\rangle_M$ and its dual $\langle\text{-},\text{-}\rangle_M$. Using this notation the dual metric $\langle\text{-},\text{-}\rangle_M$ is given explicitly by:

$$\langle \alpha, \beta \rangle_M := \begin{cases} \langle \pi^\sharp \alpha, \pi^\sharp \beta \rangle_{\mathcal{F}} & \text{if } \alpha, \beta \in (\operatorname{Ker}\pi^\sharp)^\perp, \\ \langle \alpha, \beta \rangle^\perp & \text{if } \alpha, \beta \in \operatorname{Ker}\pi^\sharp, \\ 0 & \text{in other cases.} \end{cases}$$

A simple verification, using Proposition 3.12 $(v)$ and the relation (1.1), shows that $\nabla^\bullet \pi = 0$, where $\nabla^\bullet$ comes from $\langle\text{-},\text{-}\rangle_M$. To finish the proof we construct the (smooth) contravariant tensor field $\mathcal{J}$ as the direct sum of $J$ acting along the leaves and the null map acting along the transversal directions as in the proof of Theorem 3.3 $(iii)$. With this definition the triple $(\pi, \langle\text{-},\text{-}\rangle, \mathcal{J})$ has the desired compatibility. $\square$

There is another consequence for the geometry of the foliation that arises from the transversal invariance of Proposition 3.16 $(i)$. To state it, the notions of *mean curvature field*, *characteristic form* of a foliation, and *harmonic* (or *minimal*) foliation are required. A brief account of these notions in the context of foliations for pseudo-Riemannian manifolds can be found in [25].

**Theorem 3.28.** *Let $(M, \pi)$ be an oriented Poisson manifold of dimension $n \geqslant 2$ equipped with a pseudo-Riemannian metric $\langle\text{-},\text{-}\rangle$ of signature $(r, s)$, and denote by $\mathcal{F}$ its symplectic foliation. If for every $X \in \Gamma(T\mathcal{F}^\perp)$ it holds that $\mathcal{L}_X \pi = 0$, then the leaves of the symplectic foliation are minimal submanifolds of $(M, \langle\text{-},\text{-}\rangle)$.*

*Proof.* Let us denote by $[\text{-},\text{-}]_\mathrm{s}$ the Schouten bracket of multivector fields, by $\mu \in \Omega^n(M)$ the Riemannian volume form of $(M, \langle\text{-},\text{-}\rangle)$, and by $\star \colon \Omega^\bullet(M) \to \Omega^{n-\bullet}(M)$ the Hodge star operator on forms associated to the metric $\langle\text{-},\text{-}\rangle$. Suppose that the leaves of the symplectic foliation $\mathcal{F}$ of $(M, \pi)$ are $2k$-dimensional submanifolds and denote by $\chi_{\mathcal{F}}$ the characteristic form of $\mathcal{F}$, which is precisely the Riemannian volume $2k$-form of the leaves associated to the induced metric on them. By the pseudo-Riemannian version of a result of Rummler (see [25, Cor. 3]) we have for any $X \in \Gamma(T\mathcal{F}^\perp)$ that

$$\mathcal{L}_X \chi_{\mathcal{F}} = (-1)^{s+1} \kappa(X) \chi_{\mathcal{F}} + \eta, \tag{3.9}$$

where $\kappa \in \Omega^1(M)$ denotes the mean curvature one-form of $\mathcal{F}$ and $\eta \in \Omega^{2k}(M)$ is some $2k$-form with the property: $i_{Y_1} \cdots i_{Y_{2k}} \eta = 0$ for any $Y_1, \ldots, Y_{2k} \in \Gamma(T\mathcal{F})$.

On the other hand, if $\pi^k := \pi^{\wedge k} \in \mathfrak{X}^{2k}(M)$ denotes the $k$-fold wedge product of $\pi$ with itself, we know that the rank of $\pi_p$ is $2k$, at a point of $p \in M$, if and only if $\pi_p^k \neq 0$ but $\pi_p^{k+1} = 0$. Here $2k$ is the dimension of the symplectic leaf $\mathcal{F}_p$ passing through $p$ and $\ell = n - 2k$ its codimension in $M$. Now, using the $\ell$-form $\nu := i_{\pi^k}\mu$, we can write $\mu = \chi_\mathcal{F} \wedge \nu$, where $\chi_\mathcal{F} = \star \nu$ and by the definition of the $\mu$-divergence of $X \in \Gamma(T\mathcal{F}^\perp)$ we have

$$(\operatorname{div}_\mu X)\mu = \mathcal{L}_X \mu = (\mathcal{L}_X \chi_\mathcal{F}) \wedge \nu + \chi_\mathcal{F} \wedge \mathcal{L}_X \nu. \tag{3.10}$$

From the multivector Cartan formula $i_{[X,\pi^k]_{\mathrm{S}}} = [\mathcal{L}_X, i_{\pi^k}]$ for any $X \in \mathfrak{X}(M)$, we find that

$$i_{[X,\pi^k]_{\mathrm{S}}}\mu = \mathcal{L}_X(i_{\pi^k}\mu) - (\operatorname{div}_\mu X)\, i_{\pi^k}\mu. \tag{3.11}$$

By hypothesis $\mathcal{L}_X \pi = 0$ when $X \in \Gamma(T\mathcal{F}^\perp)$, hence $[X, \pi^k]_{\mathrm{S}} = \mathcal{L}_X \pi^k = 0$, and we conclude that $\mathcal{L}_X \nu = (\operatorname{div}_\mu X)\nu$. From (3.10) this implies that $(\mathcal{L}_X \chi_\mathcal{F}) \wedge \nu = 0$.

Since $i_Y \nu = 0$ for any $Y \in T\mathcal{F}$, for any local orthonormal frame of $M$ divided in $\mathcal{F}$-tangent $E_1, \ldots, E_{2k} \in \Gamma(T\mathcal{F})$ and $\mathcal{F}$-normal $E_{2k+1}, \ldots, E_{2k+\ell} \in \Gamma(T\mathcal{F}^\perp)$ vector fields, we have from (3.9) and the properties of $\nu$, $\eta$, and $\mu = \chi_\mathcal{F} \wedge \nu$ that:

$$\begin{aligned} 0 &= \big((-1)^{s+1}\kappa(X)\,\chi_\mathcal{F} \wedge \nu + \eta \wedge \nu\big)(E_1, \ldots, E_{2k}, E_{2k+1}, \ldots, E_{2k+\ell}) \\ &= (-1)^{s+1}\kappa(X)\,\mu(E_1, \ldots, E_n) \\ &\quad + \sum_{\sigma \in S(2k,\ell)} \operatorname{sgn}(\sigma)\,\eta(E_{\sigma(1)}, \ldots, E_{\sigma(2k)})\,\nu(E_{\sigma(2k+1)}, \ldots E_{\sigma(2k+\ell)}) \\ &= (-1)^{s+1}\kappa(X). \end{aligned}$$

Since this holds for any vector field $X$ orthogonal to the leaves, we obtain that $\kappa = 0$ and by definition the leaves are minimal submanifolds of $M$. $\square$

As an immediate consequence, manifolds equipped with Riemannian and transversally invariant Poisson structures, and in particular all Kähler–Poisson or Killing–Poisson manifolds, carry *harmonic* symplectic foliations, in the sense of minimality of the leaves. Such foliations are also called *geometrically taut* (e.g. in [12]) or *minimal* although this last terminology is in conflict with other notions of minimality in foliation theory like that of a foliation with dense leaves.

## 4. Riemannian submersions and compatibilities

In this section we study the behavior of Riemann–Poisson, Kähler–Poisson, and Killing–Poisson manifolds under structure preserving surjective submersions $t \colon P \to M$. Throughout this section we consider Poisson and metric structures on $P$ and $M$ denoted by triples $(P, \pi_P, \langle \cdot, \cdot \rangle_P)$ and $(M, \pi_M, \langle \cdot, \cdot \rangle_M)$, respectively.

First, let us recall the notion of a Riemannian submersion.

**Definition 4.1.** A surjective submersion $t \colon (P, \langle \cdot, \cdot \rangle_P) \to (M, \langle \cdot, \cdot \rangle_M)$ between pseudo-Riemannian manifolds is called a *pseudo-Riemannian submersion* if the differential map $dt \colon (\operatorname{Ker} dt)^\perp \to TM$ is an isometry.

Some direct consequences are the following

**Lemma 4.2.** *Suppose that* $t \colon (P, \langle \cdot, \cdot \rangle_P) \to (M, \langle \cdot, \cdot \rangle_M)$ *is a Riemannian submersion,* $f \in C^\infty(M)$ *and* $\alpha, \beta \in \Gamma(T^*M)$. *Then*

(i) *The gradients* $\nabla(f \circ t)$ *in* $P$, *and* $\nabla f$ *in* $M$, *are* $t$-*related vector fields.*
(ii) $(t^*\alpha)^{\sharp_P} \in (\operatorname{Ker} dt)^\perp$ *and* $dt\,(t^*\alpha)^{\sharp_P} = \alpha^{\sharp_M}$.
(iii) $\langle \alpha, \beta \rangle_M \circ t = \langle t^*\alpha, t^*\beta \rangle_P$.

*Proof.* (*i*) For an arbitrary point $q \in M$ and some $p \in t^{-1}\{q\} \subset P$ let us define the subspaces $V_p P := T_p(t^{-1}\{q\}) = \operatorname{Ker} dt_p$ and $H_p P := (V_p P)^\perp$ of $T_p P$ hence $T_p P = H_p P \oplus V_p P$. Since $t$ is a submersion, it defines an isomorphism $H_p P \cong T_q M$. To avoid reference to the points we consider the vertical $VP := T(t^{-1}(P)) = \operatorname{Ker} dt$ and horizontal $HP := (VP)^\perp$ subbundles of $TP$ that define the $\langle -,- \rangle_P$-orthogonal decomposition $TP = VP \oplus HP$. Given a vector $w \in TM$, let us assume that $w = dt(v)$ for some $v \in TP$ which is decomposable as $v = v^{\text{V}} + v^{\text{H}}$, with $v^{\text{V}} \in VP$ and $v^{\text{H}} \in HP$. By definition $dt(v^{\text{V}}) = 0$, therefore, using that $t$ is a Riemannian submersion we can write

$$\begin{aligned}\langle \nabla f, w \rangle_M &= \langle \nabla f, dt(v) \rangle_M = \langle \nabla f, dt(v^{\text{H}}) \rangle_M = df(dt(v^{\text{H}})) = d(f \circ t)(v^{\text{H}}) \\ &= \langle \nabla(f \circ t), v^{\text{H}} \rangle_P = \langle dt(\nabla(f \circ t)), dt(v^{\text{H}}) \rangle_M \\ &= \langle dt(\nabla(f \circ t)), w \rangle_M.\end{aligned}$$

Since this holds for any $w \in TM$ at a fixed arbitrary point of $M$, using that the pseudo-Riemannian metric is non-degenerate, it follows $dt(\nabla(f \circ t)) = \nabla f$.

(*ii*) By definition of the bundle map $\sharp_P$, we know that $\langle (t^*\alpha)^{\sharp_P}, - \rangle_P = t^*\alpha$. Hence, for any $Z \in \operatorname{Ker} dt$ we have $\langle (t^*\alpha)^{\sharp_P}, Z \rangle_P = \alpha(dt(Z)) = 0$, which means that $(t^*\alpha)^{\sharp_P} \in (\operatorname{Ker} dt)^\perp$. The second equality follows from the fact that when restricted to $(\operatorname{Ker} dt)^\perp$ the bundle map $dt$ is an isometry.

(*iii*) Using the contravariant metric (2.3) and the previous property we get for each point $p \in P$ that

$$\langle \alpha, \beta \rangle_M |_{t(p)} = \langle dt\,(t^*\alpha)^{\sharp_P}, dt\,(t^*\beta)^{\sharp_P} \rangle_M |_{t(p)} = \langle t^*\alpha, t^*\beta \rangle_P |_p. \qquad \square$$

If in addition we assume that $t \colon P \to M$ is a Poisson map, we get

**Proposition 4.3.** *Let $\mathcal{J}_M \colon T^*M \to T^*M$ and $\mathcal{J}_P \colon T^*P \to T^*P$ be two bundle maps compatible[2] with the metric and Poisson structures of $(M, \pi_M, \langle -,- \rangle_M)$ and $(P, \pi_P, \langle -,- \rangle_P)$. For any $\alpha, \beta, \gamma \in \Gamma(T^*M)$ the following relations hold:*

$$\langle t^*\alpha, \mathcal{J}_P(t^*\beta) \rangle_P = \langle t^*\alpha, t^*(\mathcal{J}_M\beta) \rangle_P, \tag{4.1}$$

$$\langle \nabla^{t^*\alpha}(t^*\beta), t^*\gamma \rangle_P = t^*\langle \nabla^{\alpha}\beta, \gamma \rangle_M, \tag{4.2}$$

$$\pi_P(\nabla^{t^*\alpha}(t^*\beta), t^*\gamma) = t^*\bigl(\pi_M(\nabla^{\alpha}\beta, \gamma)\bigr). \tag{4.3}$$

*Proof.* The first identity comes from

$$\langle t^*\alpha, \mathcal{J}_P(t^*\beta) \rangle_P = \pi_P(t^*\alpha, t^*\beta) = t^*\bigl(\pi_M(\alpha, \beta)\bigr) = t^*\bigl(\langle \alpha, \mathcal{J}_M\beta \rangle_M\bigr) = \langle t^*\alpha, t^*(\mathcal{J}_M\beta) \rangle_P.$$

Relation (4.2) follows from a direct calculation by using (1.1) and the properties in Lemma 4.2. Finally, for (4.3) we have

$$\pi_P(\nabla^{t^*\alpha}(t^*\beta), t^*\gamma) = \langle t^*(\nabla^{\alpha}\beta), t^*\gamma' \rangle_P = \pi_P(t^*(\nabla^{\alpha}\beta), t^*\gamma),$$

where $t^*\gamma' = \mathcal{J}_P(t^*\gamma)$, and the conclusion follows. $\square$

Note that for one-forms $\alpha, \beta \in \Gamma(T^*M)$, the first claim of Proposition 4.3 implies that $\langle t^*\alpha, \mathcal{J}_P(t^*\beta) \rangle_P$ is a $t$-basic function and, therefore, $\mathcal{J}_P\bigl(\Gamma(t^*T^*M)\bigr) \subset \Gamma(t^*T^*M)$. Another result on $t$-basic forms is the following equivalence for partially complex structures preserved by Riemannian submersions and Poisson maps,

**Proposition 4.4.** *Let $(P, \pi_P, \langle -,- \rangle_P, \mathcal{J}_P)$ be an almost Kähler–Poisson manifold. The partially complex structure $\mathcal{J}_P$ preserves basic one-forms on $P$ if and only if there exists a partially complex structure $\mathcal{J}_M$ on $M$ so that $(M, \pi_M, \langle -,- \rangle_M, \mathcal{J}_M)$ is an almost Kähler–Poisson manifold satisfying $t^* \circ \mathcal{J}_M = \mathcal{J}_P \circ t^*$. Moreover, relations (4.2) and (4.3) hold in this case.*

---

[2] here, $\mathcal{J}$ just satisfies $\pi^\sharp + \sharp \circ \mathcal{J} = 0$ but may not be a contravariant $f$-structure.

*Proof.* First we will assume that $J_M$ exists. In this case, for any one-form $\alpha$ in $M$ we can verify that $J_P(t^*\alpha) = t^*(J_M\alpha)$, that is $(J_P \circ t^*)(\alpha)$ is again basic. For the converse, note that to any $\alpha \in \Omega^1(M)$ corresponds a unique $\alpha' \in \Omega^1(M)$ so that $t^*\alpha' = J_P(t^*\alpha)$. The bundle map $J_M \colon T^*M \to T^*M$ is defined via the relation $J_M\alpha = \alpha'$, i.e.,

$$t^* \circ J_M = J_P \circ t^*. \tag{4.4}$$

From the definition we can verify that $t^* \circ (J_M^3 + J_M) = (J_P^3 + J_P) \circ t^* = 0$, hence $J_M$ is partially complex. Using the compatibilities of the structures in $P$ and the fact that $t$ is a Poisson map and a Riemannian submersion we get

$$t^*\pi_M(\alpha, \beta) = \langle t^*\alpha, J_P(t^*\beta)\rangle_P = \langle t^*\alpha, t^*(J_M\beta)\rangle_P = t^*\langle \alpha, J_M\beta\rangle_M,$$

hence $\pi_M(\alpha, \beta) = \langle \alpha, J_M\beta\rangle_M$, which means that $(M, \pi_M, \langle \text{-},\text{-}\rangle_M, J_M)$ is almost Kähler–Poisson. Finally, the relation $\langle t^*\alpha, J_P(t^*\beta)\rangle_P = \langle t^*\alpha, t^*(J_M\beta)\rangle_P$ is trivially satisfied, and as a direct consequence of the previous proof we obtain that (4.2) and (4.3) hold. □

4.1. **Killing–Poisson case.** For a Poisson map $t \colon (P, \pi_P) \to (M, \pi_M)$ we study sufficient conditions on the structures $\pi_P$ and $\pi_M$ so that the properties given in Proposition 3.16 are preserved by $t$. First, we consider the transversal invariance of the Poisson bivector fields or, equivalently, that gradient vector fields associated to Casimir functions are Poisson vector fields.

**Proposition 4.5.** *Let $t \colon P \to M$ be a Riemannian submersion and a Poisson map simultaneously and suppose that $Z(M, \pi_M)$ and $Z(P, \pi_P)$ are the centers of $\Omega^1(P)$ and $\Omega^1(M)$, respectively, regarded as Lie algebras with their Koszul brackets. If $(P, \pi_P, \langle \text{-},\text{-}\rangle_P)$ is transversally invariant as in Proposition 3.16 (i) and one of the following conditions hold:*

*(i)* $t^*\big(Z(M, \pi_M)\big) \subset Z(P, \pi_P)$,
*(ii)* $t^*(\operatorname{Ker} \pi_M^\sharp) \subset \operatorname{Ker} \pi_P^\sharp$,

*then $(M, \pi_M, \langle \text{-},\text{-}\rangle_M)$ also satisfies the transversal invariance property.*

*Proof.* For condition (i) we employ the equivalences stated in [6, Prop 2.4] and the fact that Proposition 3.16 (i) holds if and only if $\nabla^\bullet \alpha = 0$ for any $\alpha \in \Omega^1(M)$ that satisfies (4.2) and belongs to the center $Z(M, \pi_M)$ of the Lie algebra $\big(\Omega^1(M), [\text{-},\text{-}]_{\pi_M}\big)$. Clearly, the conclusion also holds for condition (ii) assuming that

$$t_*(\mathcal{L}_X \pi_P) = \mathcal{L}_Y \pi_M$$

for $t$-related vector fields $X \in \Gamma(TP)$ and $Y \in \Gamma(TM)$. □

*Example* 4.6. When $M$ and $P$ are symplectic manifolds and $t \colon P \to M$ is both, a Riemannian submersion and a symplectic map, the Killing–Poisson condition on $P$ implies the Killing–Poisson condition on $M$. (For condition (ii) see Example 3.19.) ◇

*Example* 4.7. If $(M, \pi_M)$ and $(N, \pi_N)$ are Poisson manifolds, the product Poisson structure on $P = M \times N$ satisfies condition (ii) of Proposition 4.5. Moreover, the product metric on $P$ turns $t \colon P \to M$ into a Riemannian submersion, hence we get the result. ◇

4.2. **Riemann–Poisson and Kähler–Poisson case.** The main result of this section describes the behavior of Kähler–Poisson structures under Riemannian submersions. An important application is to the reduction of Kähler manifolds under a group of symmetries. In particular, we will study the case of compact Lie groups. For the remaining of this section we consider two smooth manifolds carrying Poisson and Riemannian structures $(M, \pi_M, \langle \text{-},\text{-}\rangle_M)$ and $(P, \pi_P, \langle \text{-},\text{-}\rangle_P)$ with corresponding cometrics $\langle \text{-},\text{-}\rangle_M$ and $\langle \text{-},\text{-}\rangle_P$.

**Theorem 4.8.** *Let $t: P \to M$ be a Riemannian submersion that is also a Poisson map, then $\nabla^\bullet \pi_M = 0$ whenever $\nabla^\bullet \pi_P = 0$. Furthermore, if there exist partially complex structures $J_M$ on $M$ and $J_P$ on $P$ such that the triples $(\pi_M, \langle -,-\rangle_M, J_M)$ and $(\pi_P, \langle -,-\rangle_P, J_P)$ are themselves compatible then $M$ is Kähler–Poisson whenever $P$ is Kähler–Poisson.*

*Proof.* For the first claim just note that from (4.2), (4.3) and using that $t^*$ is injective it follows that
$$t^*(\pi_M^\sharp(\alpha)\,\pi_M(\beta,\gamma)) = t^*(\pi_M(\nabla^\alpha \beta, \gamma) + \pi_M(\beta, \nabla^\alpha \gamma))$$
for all one-forms $\alpha, \beta, \gamma \in \Gamma(T^*M)$. Finally, for the second claim it is enough to use the equivalence in Proposition 3.7. □

The converse is not true, as we will see in the case of a non-flat Riemannian manifold. Let $M$ be such a manifold equipped with the null Poisson structure $\pi = 0$, which trivially satisfies (3.3). If we fix the manifold $P := T^*M$ with its canonical symplectic form $\omega_{\mathrm{can}}$, then the natural projection $t: T^*M \to M$ is a Poisson map but there does not exist a metric on $P$ turning it into a Kähler manifold [28].

*Example* 4.9. Suppose that a symplectic groupoid $(\mathcal{G}, \omega) \rightrightarrows M$ is equipped with a metric $\langle -,-\rangle_\mathcal{G}$ so that $(\mathcal{G}, \omega, \langle -,-\rangle_\mathcal{G})$ is Kähler. If in addition we suppose that $\langle -,-\rangle_\mathcal{G}$ is invariant under left translations $L_g: t^{-1}(s(g)) \to t^{-1}(t(g))$, the Lie groupoid will be called a *Kähler groupoid*. With this assumption we obtain a metric on $M$ for which the target map $t$ is a Riemannian submersion (the same argument holds for the source map but assuming that the metric is right-invariant). As a direct consequence of Proposition 4.8 we get that $(M, \pi)$ is Kähler–Poisson just from the fact that $t$ is a Poisson map (or that $s$ is an anti-Poisson map). ◇

*Example* 4.10 (Example 3.6 revisited). Consider $P = \mathbb{R}^4$ with its canonical Kähler structure $(\widetilde{\omega}, J_P, \langle -,-\rangle_P)$. Let $M = \mathbb{R}^3$ be endowed with the Poisson bivector $\pi_M = \partial_1 \wedge \partial_2$ and denote by $t: \mathbb{R}^4 \to \mathbb{R}^3$ the projection onto the first three components. Since $t$ is a Poisson map and a Riemannian submersion with the aid of Theorem 4.8 we conclude that $M$ is a Riemann–Poisson manifold. Furthermore, $\mathbb{R}^3$ is a Kähler–Poisson manifold if we consider the compatible partially complex structure $J_M = \partial_1 \otimes dx_2 - \partial_2 \otimes dx_1$ as in Example 3.6. Note that the partially complex $J_P$ does not come as a $t$-pull-back (4.5) of $J_M$ because $J_P \circ dt^*$ does not produce basic one-forms (see Proposition **??**). Nevertheless, $J_M$ can be obtained from $J_P$ by the formula
$$J_M := -\flat_M \circ dt \circ J_P \circ \sharp_P \circ dt^*, \tag{4.5}$$
where $\flat_M$ denotes the flat isomorphism associated to the metric $\langle -,-\rangle_M$ of $M$. ◇

The previous constructions can be extended to *cosymplectic manifolds*, of which $\mathbb{R}^3$ is the basic example. A cosymplectic manifold is a triple $(M, \omega, \eta)$ where $M$ is a $2n+1$-dimensional manifold equipped with a closed two-form $\omega$ and a closed one-form $\eta$ for which $\omega^n \wedge \eta \neq 0$. In this case, $P = M \times \mathbb{R}$ is symplectic manifold with symplectic form $\widetilde{\omega} = t^*\omega + t^*\eta \wedge ds$ where $t: P \to M$ is the usual projection map. Moreover, $M$ can be equipped with a Poisson bivector $\pi_M$ for which $t$ is a Poisson map (for details see [14]). If $\langle -,-\rangle_M$ is a metric on $M$ whose dual is $\langle -,-\rangle_M$, and $\langle -,-\rangle_P = \langle -,-\rangle_M + \langle -,-\rangle_\mathbb{R}$ is the metric on $P$, with dual metric $\langle -,-\rangle_P = \langle -,-\rangle_M + \langle -,-\rangle_\mathbb{R}$, we get that $t$ is also a Riemannian submersion. A direct computation using the Levi-Civita connection on $P$, and that the one-form $ds$ in $\mathbb{R}$ is parallel, yields the following corollary of Theorem 4.8:

**Corollary 4.11.** *If $\omega$ and $\eta$ are parallel, then $(P, \widetilde{\omega}^{-1})$ and $(M, \pi_M)$ are Riemann–Poisson.*

What remains to verify is the existence of an almost partially complex structure compatible with these metrics and Poisson structures in the same way as it was done in the previous example for the cosymplectic manifold $\mathbb{R}^3$.

**Proposition 4.12.** *Consider the same situation of the previous corollary and also assume that the metric on $M$ is normalized so that $\|\eta\| = 1$. If $J_M := \sharp_M \circ \omega^\sharp$ satisfies the conditions*

(i) $\mathrm{Im}(J_M) \subset \mathrm{Ker}(\eta)$,
(ii) $J_M(\sharp_M \eta) = 0$,
(iii) $(J_M^2 + \mathrm{Id})(X) = \eta(X)\sharp_M \eta$, *for all $X \in \Gamma(TM)$,*

*then $J_P = \sharp_P \circ \widetilde{\omega}^\sharp$ is a complex structure, $\mathcal{J}_M$ defined as in (4.5) is a partially complex structure, $(P, \widetilde{\omega}, J_P, \langle\text{-},\text{-}\rangle_P)$ is a Kähler manifold, and $(M, \pi_M, \mathcal{J}_M, \langle\text{-},\text{-}\rangle_M)$ is Kähler–Poisson.*

*Proof.* Conditions (i)–(iii) imply that $J_P(X + a\partial_s) = J_M(X) - a(\sharp_M \eta) + \eta(X)\,\partial_s$ is a compatible complex structure, and the previous corollary together with Proposition 3.7 lead to the conclusion that $(P, \widetilde{\omega}, J_P, \langle\text{-},\text{-}\rangle_P)$ is a Kähler manifold. For the second claim we must prove that $\mathcal{J}_M$ is compatible and partially complex. To prove compatibility we use the equalities $-\flat_P \circ J_P \circ \sharp_P = J_P = -\flat_P \circ \pi_P^\sharp$ to verify that $\sharp_M \circ \mathcal{J}_M = -\pi_M^\sharp$. To prove the relation $\mathcal{J}_M^3 + \mathcal{J}_M = 0$, recall that at each point $m \in M$ we have the splitting $T_m^* M = \mathrm{Im}(\omega_m^\sharp) \oplus \langle \eta_m \rangle$, thus by condition (ii) we get that $\mathcal{J}_M(\eta) = 0$, and from the fact that $t$ is a Riemannian submersion $\mathcal{J}_M^2 + \mathrm{Id} = 0$ when restricted to $\mathrm{Im}(\omega_m^\sharp)$. Finally, $(M, \pi_M, \mathcal{J}_M, \langle\text{-},\text{-}\rangle_M)$ is Kähler–Poisson as consequence of the previous corollary and Proposition 3.7. $\square$

## 5. Riemann–Poisson and Kähler–Poisson symmetries

In this section we will assume that there is a smooth action $\varphi : G \times P \to P$ on a Poisson and Riemannian manifold $(P, \pi_P, \langle\text{-},\text{-}\rangle_P)$ in such a way that each diffeomorphism $\varphi_g$ is a Poisson map and an isometry. In the case of proper and free actions, these two assumptions yield that the quotient map $M := P/G$ inherits a unique Poisson structure $\pi_M$ and a Riemannian metric $\langle\text{-},\text{-}\rangle_M$, so that the quotient map $t: P \to M$ is a Poisson map and a Riemannian submersion.

In the case of a Kähler–Poisson manifold $(P, \pi_P, \langle\text{-},\text{-}\rangle_P, \mathcal{J}_P)$, the preserving conditions on $\pi_P$ and $\langle\text{-},\text{-}\rangle_P$ imply that $d\varphi_g^* \circ \mathcal{J}_P = \mathcal{J}_P \circ d\varphi_g^*$, for all $g \in G$. It is routine to verify that such a commuting relation yields[3] $\mathcal{L}_{u_P}(\mathcal{J}_P(t^*\alpha)) = 0$ for all one-forms $\alpha$ on $M$, and any $u \in \mathfrak{g} = \mathrm{Lie}(G)$. Hence, the natural question that arises now is under which conditions we get that the basic one-forms $\Omega_b^1(P)$ on $P$ are preserved by $\mathcal{J}_P$.

**Lemma 5.1.** *The partially complex structure $\mathcal{J}_P$ preserves basic one-forms in $\Omega_b^1(P)$ if and only if $J_P \mathcal{V} \subset \mathcal{V}$, where $J_P := -\sharp \circ \mathcal{J}_P \circ \flat$, and $\mathcal{V}$ is the vertical space of the quotient map $t$.*

*Proof.* Recall that the one-form $\gamma$ in $P$ is basic with respect to $t$ if and only if $\mathcal{L}_{u_P}\gamma = 0$ and $i_{u_P}\gamma = 0$. Hence, by the comment made before the statement of the lemma, we note that $\mathcal{J}_P \Omega_b^1(P) \subset \Omega_b^1(P)$ if and only if $\mathcal{J}_P' \mathcal{V} \subset \mathcal{V}$ where $\mathcal{J}_P' : TP \to TP$ is the dual bundle map of $\mathcal{J}_P$. Thus, the claim is proven if we verify that $\mathcal{J}_P' = J_P$.

Note that, by definition, we have the following facts

$$\flat' = \flat \quad \text{and} \quad \mathcal{J}_P' \circ \sharp' = \pi_P^\sharp,$$

and both yield the relation $\sharp \circ \mathcal{J}_P' = -\pi_P^\sharp$. Computing this relation for all $\alpha, \beta \in T^*P$ we get

$$\langle \beta, \mathcal{J}_P \alpha \rangle = \beta(-\pi_P^\sharp \alpha) = \beta(\sharp \circ \mathcal{J}_P') = \langle \beta, \mathcal{J}_P' \alpha \rangle$$

which finally says that $\mathcal{J}_P' = J_P$ and the result is proved. $\square$

The previous situation allows us to apply Theorem 4.8 in order to obtain

**Theorem 5.2.** *Under the presence of a $G$-action with the assumptions stated above:*

---

[3] we denote by $u_P$ the infinitesimal generator of $u \in \mathfrak{g}$ by the $G$-action on $P$ and recall that the space vector bundle spanned by $u_P$ coincides with the vertical space of the quotient map $t$.

(i) If $(P, \pi_P, \langle-,-\rangle_P)$ is a Riemann–Poisson manifold, then the reduced manifold $P/G$ is again a Riemann–Poisson manifold.

(ii) If $(P, \pi_P, \langle-,-\rangle_P, J_P)$ is a Kähler–Poisson manifold and $J_P \mathcal{V} \subset \mathcal{V}$, then the reduced manifold $P/G$ is again a Kähler–Poisson manifold.

As a direct consequence of (i), we get that the quotient of a Kähler manifold by the action of a group of preserving symmetries and symplectomorphisms is a Riemann–Poisson manifold. Therefore, each leaf is again a Kähler manifold but, nevertheless, the metric is not necessarily the restricted one (see Theorem 3.24).

Another interesting consequence of the previous theorem is related to the Kähler reduction by a Hamiltonian $G$-action on a regular value $\zeta \neq 0$ of the moment map $\mu : P \to \mathfrak{g}^*$ (see the comment after the proof of [10, Ch. 8, Thm. 3]). In the discussion, the sufficient condition in order to guarantee the existence of a Kähler structure by the Marsden–Weinstein reduction on $\mu^{-1}(\zeta)/G_\zeta$ is that the bundle $\mathcal{V}_\zeta^\perp$ must be invariant under $J_\zeta$, the restriction of $J_P$ to the submanifold $\mu^{-1}(\zeta)$, where $\mathcal{V}_\zeta$ is the vertical space of the action of the isotropy group $G_\zeta$ on $\mu^{-1}(\zeta)$. As a consequence of the hypothesis $J_P \mathcal{V} \subset \mathcal{V}$ we can verify that $J_\zeta \mathcal{V}_\zeta \subset \mathcal{V}_\zeta$, and from this it follows that $\mathcal{V}_\zeta^\perp$ is also $J_\zeta$-invariant. Thus, the hypothesis $J_P \mathcal{V} \subset \mathcal{V}$ ensures that Kähler reduction works on regular values of the moment map for a Hamiltonian action on a Kähler manifold $P$.

*Example* 5.3. Going back to Example 3.6, for integers $n > 2$ and $r, s \in \{1, \ldots, n\}$ with $r < s$, let us denote here by $\mathbb{R}^n_{(rs)} := (\mathbb{R}^n, \pi_{(rs)}, \langle-,-\rangle, J_{(rs)})$ the Kähler–Poisson structure defined there. Taking $m \leqslant n - 2$ and increasingly ordered indices $i_1, \ldots, i_m \in \{1, \ldots, n\} \setminus \{r, s\}$, we define the space $\mathcal{V}_m := \operatorname{span}\{\partial_{i_1}, \ldots, \partial_{i_m}\} \cong \mathbb{R}^m$ and consider the (free and proper) $\mathbb{R}^m$-action on $\mathbb{R}^n$ by translations: $\tau : \mathbb{R}^n \times \mathbb{R}^m \to \mathbb{R}^n$, $(p, v) \mapsto \tau_v(p) := p + v$. Note that this action is an isometry for the Euclidean cometric and also preserves $\pi_{(rs)}$, hence commutes with $J_{(rs)}$. The quotient map $t : \mathbb{R}^n_{(rs)} \to \mathbb{R}^{n-m}_{(rs)} := \mathbb{R}^n_{(rs)}/\sim_\tau$, is the projection onto the remaining coordinates $(x_{j_1}, \ldots, x_{j_m})$, where $j_1, \ldots, j_m \in \{1, \ldots, n\} \setminus \{i_1, \ldots, i_m\}$. Moreover the vertical space is precisely $\mathcal{V}_m$, which is $J'_{(rs)}$-invariant, $J'_{(rs)} := -(dx^r \otimes \partial_s - dx^s \otimes \partial_r)$. This discussion allows us to apply Theorem 5.2, and obtain that the quotient $\mathbb{R}^{n-m}_{(rs)}$ is a Kähler–Poisson manifold, as expected. In the case $n - m = 2$, this quotient is precisely the Kähler manifold $\mathbb{C} \cong \mathbb{R}^2$ and, in particular for $n = 3$, $m = 1$ we have the foliation of $\mathbb{R}^3$ by $\mathbb{C}$-planes as discussed at the end of Example 3.6. $\diamondsuit$

In order to present another relevant instance of Kähler–Poisson reduction we consider the following situation: let $(P, \omega, \langle-,-\rangle, J)$ denote a Kähler manifold with a Hamiltonian $G$-action preserving these structures, carrying an associated moment map $\mu$, and where all the actions involved are free and proper.

**Proposition 5.4.** *If under the previous geometrical setting we assume, in addition, that each $\zeta \in \operatorname{Im}(\mu)$ is a regular value for $\mu$ and $H_\zeta := (\mathcal{V} \cap T\mu^{-1}(\zeta))^\perp$ is $J|_{T\mu^{-1}(\zeta)}$-invariant, then the orbit space $M = P/G$ is Kähler–Poisson.*

*Proof.* The condition on $H_\zeta$ implies that the Marsden–Weinstein symplectic reduction $(P_\zeta := \mu^{-1}(\zeta)/G_\zeta, \omega_\zeta)$ is a Kähler manifold (see [10, Ch. 8]) with compatible complex structure denoted by $J_\zeta$. From the compatibilities on each reduced manifold we have that:

$$J_\zeta = \sharp_\zeta \circ \omega_\zeta^\flat, \quad \text{and} \quad J_\zeta \circ \sharp_\zeta = -\sharp_\zeta \circ J_\zeta.$$

We also know (from Theorem 5.2) that $M$ is Riemann–Poisson with Poisson bivector $\pi_M$ and on each leaf of the symplectic foliation $\mathcal{F}$ the symplectic form coincides with $\omega_\zeta$, and the metric $\langle-,-\rangle_M$ restricts to $P_\zeta$, which is a symplectic manifold sitting inside $M$ as a union of symplectic leaves of the reduced Poisson manifold $M$.

On $M$ it is possible to construct a smooth bundle map $\mathcal{J}_M : T^*M \to T^*M$ as the composition $\mathcal{J}_M := -\flat_M \circ \pi_M^\sharp$, thus $(\pi_M, \langle\text{-},\text{-}\rangle_M, \mathcal{J}_M)$ is a compatible triple. Furthermore, when we restrict to $T^*\mathcal{F}$ we can verify that

$$\mathcal{J}_M|_{T^*\mathcal{F}} = -\flat|_{T\mathcal{F}} \circ \pi_M^\sharp|_{T^*\mathcal{F}} = \flat_\zeta \circ (\omega_\zeta^\flat)^{-1}$$
$$= -\flat_\zeta \circ J_\zeta \circ \sharp_\zeta = \flat_\zeta \circ \sharp_\zeta \circ \mathcal{J}_\zeta = \mathcal{J}_\zeta,$$

which means that $\mathcal{J}_M$ coincides with the complex structure $\mathcal{J}_\zeta$ on each leaf of $\mathcal{F}$. From this claim and the fact that $\text{Ker}(\mathcal{J}_M) = \text{Ker}(\pi_M^\sharp)$, we conclude that $\mathcal{J}_M^3 + \mathcal{J}_M = 0$, i.e., $(M, \pi_M, \langle\text{-},\text{-}\rangle_M, \mathcal{J}_M)$ is an almost Kähler–Poisson manifold. Finally, as $M$ is Riemann–Poisson, from $\nabla^\bullet \pi = 0$ we get that $(M, \pi_M, \langle\text{-},\text{-}\rangle_M, \mathcal{J}_M)$ is indeed a Kähler–Poisson manifold (see Proposition 3.7). □

The case of the group $G = S^1$ acting on $P = \mathbb{C}_0^n$ does not fit into Theorem 5.2 because of the rank of the vertical space $\mathcal{V}$. But since $G$ is abelian we are in the situation of the previous proposition (again, see the comment after the proof of [10, Ch. 8, Thm. 3]) which leads us to conclude that $\mathbb{C}_0^n/S^1$ is Kähler–Poisson.

5.1. **Compact Lie groups.** Let $G$ be a compact Lie group. The aim here is to study conditions for the existence of a Riemann–Poisson structure associated to the linear Poisson structure in the dual $\mathfrak{g}^*$ of the Lie algebra $\mathfrak{g} := \text{Lie}(G)$. We denote by $B$ the Killing form of $\mathfrak{g}$, which is a non-degenerate, Ad-invariant symmetric bilinear form[4]. As it is done in [3, 4, 18], the manifold $P = G \times \mathfrak{g}^* \cong T^*G$ is equipped with a Kähler structure (i.e., it is Kähler–Poisson) with an underlying symplectic form $\omega$ symplectomorphic to the canonical one $\omega_{\text{can}}$ in $T^*G$ via the trivialization by left translations.

In order to write the metric explicitly, we must start considering the complex structure $J$ on $TG \cong G^\mathbb{C}$, where $G^\mathbb{C} := \exp\{\mathfrak{g} + i\mathfrak{g}\}$ denotes the complexification of $G$. The complex structure on the unit $J_{(e,e)} : \mathfrak{g}^2 \to \mathfrak{g}^2$ is defined by $J_{(e,e)}(u,v) = (-v, u)$ and extending it by translation to $TG$. In this case, for any $X \in TP$ there exists a unique pair of elements $u(X), v(X) \in \mathfrak{g}$ for which

$$X = u(X)_P + Jv(X)_P.$$

The metric is then given by (see [3, 18])

$$\langle X, Y \rangle = \psi(u(X))(u(Y)) + \mu[v(X), u(Y)] + B(v(X), d^u(\mu Y)), \tag{5.1}$$

where $\mu : G \times \mathfrak{g}^* \to \mathfrak{g}^*$ is the anti-projection, $d^u := d + [u(\text{-}), \text{-}]$, and $\psi : P \to \text{Hom}(\mathfrak{g}, \mathfrak{g}^*)$ satisfies $\psi(u)(v) = \omega(u_P, Jv_P)$.

Following the constructions in [3, 18] we have that $G$ acts by isometries, therefore we have Kähler symmetries. Since the action on $P$ is trivial in the second component, the reduction is just $\mathfrak{g}^*$ with its linear Poisson structure. Moreover, $\mathfrak{g}^*$ inherits a quotient metric $\langle\text{-},\text{-}\rangle_{\text{red}}$, coming from (5.1), that turns the projection map $t$ into a Riemannian submersion (see also [3, Thm. 4.1]). In addition, note that the previous description fits into the situation of Theorem 5.2. From the construction we obtain that the vertical space of the quotient map is $\mathcal{V}|_{(e,e)} = \mathfrak{g} \times \{0\}$, and by definition of $J$ it is possible to verify that $J\mathcal{V}|_{(e,e)}$ is not contained in $\mathcal{V}|_{(e,e)}$, thus we are in the case $(i)$ of Theorem 5.2, i.e., $\mathfrak{g}^*$ is a Riemann–Poisson manifold. As consequence, each leaf inherits a Kähler structure but the metric is not the restricted one as it is observed in the literature.

*Remark* 5.5. Here, we want to compare our compatibility conditions with the case of groupoid 2-metrics, a stronger notion of metric compatibility on Lie groupoids introduced in [16]. First note that via [16, Thm. 4.3.4] we obtain that $\mathcal{G}$ has a 2-metric (because for an action groupoid $\mathcal{G} := G \ltimes M$ with compact $G$ the map $t \times s \colon \mathcal{G} \times \mathcal{G} \to M \times M$ is proper), in particular the inversion map of the groupoid $\mathrm{inv} \colon \mathcal{G} \to \mathcal{G}$ is an isometry. If such metric realizes also the Kähler structure on $\mathcal{G}$ then $d\,\mathrm{inv} \circ J + J \circ d\,\mathrm{inv} = 0$, because the inversion in the groupoid $\mathcal{G}$ is an anti-symplectomorphism. But using the fact that the complex structure at the identity is given by $J(u, v) = (-v, u)$ and computing the differential $d\,\mathrm{inv}$ of the inversion map $\mathrm{inv}$, we can verify that $d\,\mathrm{inv} \circ J + J \circ d\,\mathrm{inv} \neq 0$, which implies that the Riemannian metric on $\mathcal{G}$ given by (5.1) does not come from a 2-metric.

---

[4]Henceforth, we identify $\mathfrak{g}$ and $\mathfrak{g}^*$ under this bilinear form, referring to the appropriate one according to the context.

Nicolás Martínez Alba
Departamento de Matemáticas
Universidad Nacional de Colombia
Cra. 30 N. 45-03, Bogotá, Colombia
*E-mail:* nmartineza@unal.edu.co

Andrés Vargas
Departamento de Matemáticas
Pontificia Universidad Javeriana
Cra. 7 N. 40-62, Bogotá, Colombia
*E-mail:* a.vargasd@javeriana.edu.co